\documentclass[10pt, a4paper]{amsart}
\usepackage{amscd,amssymb,eucal,eufrak,mathrsfs,amsmath}
\usepackage{hyperref} 

\input xy
\xyoption{all}
\usepackage{epsfig}

\newtheorem{thm}{Theorem}[section]
\newtheorem{cor}[thm]{Corollary}
\newtheorem{lem}[thm]{Lemma}
\newtheorem{prop}[thm]{Proposition}

\newtheorem{question}[thm]{Question}
\theoremstyle{definition}
\newtheorem{defn}[thm]{Definition}

\theoremstyle{remark}
\newtheorem{remark}[thm]{Remark}
\newtheorem{remarks}[thm]{Remarks}
\newtheorem{example}[thm]{Example}

\numberwithin{equation}{section}

\newcommand{\delete}[1]{} 
\newcommand{\nt}{\noindent}

\def\eps{{\varepsilon}}
\def\a{\alpha}
\def\om{\omega}

\def\t{\tau}
\def\s{\sigma}
\def\g{\gamma}

\newcommand{\sk}{\vskip 0.2cm}
\newcommand{\ssk}{\vskip 0.1cm}

\newcommand{\hsk}{\hskip 0.1cm}

\newcommand{\nl}{\newline}
\newcommand{\ben}{\begin{enumerate}}

\newcommand{\een}{\end{enumerate}}
\newcommand{\bit}{\begin{itemize}}

\newcommand{\eit}{\end{itemize}}

\def\R {{\mathbb R}}
\def\N {{\mathbb N}}
\def\Z {{\mathbb Z}}

\def\T {{\mathbb T}}

\def\norm#1{\left\Vert#1\right\Vert}

\newcommand{\Acal}{\mathcal{A}}

\def\lca {locally compact abelian }
\def\nbd {neighborhood }
\newcommand{\TFAE}{The following are equivalent: }
\def\wrt{with respect to }

\def\QED{\nobreak\quad\ifmmode\roman{Q.E.D.}\else{\rm Q.E.D.}\fi}

\begin{document}

\title[]
{Every topological group is a group retract of a minimal group}

\author[]{Michael Megrelishvili}
\address{Department of Mathematics,
Bar-Ilan University, 52900 Ramat-Gan, Israel}
\email{megereli@math.biu.ac.il}
\urladdr{http://www.math.biu.ac.il/$^\sim$megereli}


\date{October 6, 2006}








\begin{abstract}
We show that every Hausdorff topological group is a group retract
of a minimal topological group. This first was conjectured by
Pestov in 1983. Our main result leads to a solution of some
problems of Arhangel'skii. One of them is the problem about
representability of a group as a quotient of a minimal group
(Problem 519 in the first edition of "Open Problems in Topology").
Our approach is based on \emph{generalized Heisenberg groups} and
on groups arising from group representations on Banach spaces and
in bilinear mappings.
\end{abstract}

\maketitle

 \tableofcontents

\setcounter{tocdepth}{1}

\section{Introduction}
\label{s:intro}

\sk

 A Hausdorff topological group $G$ is {\it minimal} (introduced by
Stephenson \cite{St} and Do\"\i chinov \cite{Do}) if $G$ does not
admit a strictly coarser Hausdorff group topology. \emph{Totally
minimal} groups are defined by Dikranjan and Prodanov \cite{DP1}
as those Hausdorff groups $G$ such that all Hausdorff quotients
are minimal (later these groups were studied also by Schwanengel
\cite{Sw} under the name {\em $q$-minimal groups}).
 First we recall some facts about minimality; mainly concerning
 the purposes of the present work. For more comprehensive
 information about minimal groups theory we refer to the book \cite{DPS}, review papers
 \cite{CHR} and \cite{Di} and also a recent article by Dikranjan and the present author \cite{DM}.

%

Unless explicitly stated otherwise, all spaces in this paper are
at least Hausdorff.
Most obvious examples of minimal groups are compact groups.
Stephenson showed \cite{St} that every abelian locally compact
minimal group must be compact. Prodanov and Stoyanov established
one of the most fundamental results in the theory proving that
every abelian minimal group is precompact \cite[Section 2.7]{DPS}.
Dierolf and Schwanengel \cite{DS} using semidirect products found
some interesting examples of non-precompact (hence, non-abelian)
minimal groups. These results imply for instance that an arbitrary
discrete group is a group retract \footnote{\emph{Group retract}
means that the corresponding retraction is a group homomorphism.}
of a locally compact minimal group; also the semidirect product
$\R \leftthreetimes \R_+$ of the group of all reals $\R$ with the
multiplicative group $\R_+$ of positive reals is minimal (now it
is known \cite{me-Gmin} that $\R^n \leftthreetimes \R_+$ is
minimal for every $n \in \N$). It follows that many minimal groups
may have nonminimal quotients (in other words, quite often minimal
groups fail to be totally minimal) and nonminimal closed
subgroups. Note that all \emph{closed} subgroups of an abelian
minimal group are again minimal
(see \cite[Proposition 2.5.7]{DPS}). Motivated by these results
Arhangel'skii posed the following two natural questions.

\ssk

\nt {\textbf{Question A:} (See \cite[Problem VI.6]{Arh}) \
\emph{Is every topological group a quotient of a minimal group~?}

\ssk

\nt {\textbf{Question B:} (See \cite[Section 3.3F, Question
3.3.1(a)]{CHR} and \cite[page 57]{Di}) \ \emph{Is every
topological group $G$ a closed subgroup of a minimal group~ $M$ ?}

\ssk

Question A appears also in the volume of "Open Problems in
Topology" \cite[Problem 519]{OPIT1} and in two review papers
\cite[Question 3.3.1]{CHR} and \cite[Question 2.9]{Di}.

\ssk The following was conjectured by Pestov.

\ssk

\nt {\textbf{Conjecture}: (Pestov 1983) \ \emph{Every topological
group is a group retract of a minimal topological group.}

\ssk

Remus and Stoyanov \cite{RS} proved that every compactly generated
locally compact abelian group is a group retract of a minimal
locally compact group. In \cite{me-min} we show that
\emph{Heisenberg type groups}
frequently are minimal (see Section \ref{s:GH}). For instance, if
$G$ is locally compact abelian with the canonical duality mapping
$\om: G^* \times G \to \T$ then the corresponding generalized
Heisenberg group $H(\om):=(\T \times G^*) \leftthreetimes G$ is
minimal. It follows that every abelian
locally compact group is a group retract of a minimal locally
compact group, as $G$ is obviously a natural retract of $H(\om)$
(see Section \ref{s:actions}).

By \cite[Theorem 4.13]{me-min} every abelian topological group is
a quotient of a minimal group. Note also that by \cite[Theorem
6.12]{me-fr} every topological subgroup $G$ of the group $Iso(V)$
of all linear isometries of a reflexive (even Asplund) Banach
space $V$ is a group retract of a minimal group. This includes all
locally compact groups because they are subgroups of $Iso(H)$
where $H$ is a Hilbert space.

In the present paper we obtain the following main result (see
Theorem \ref{t:main} below).

\vskip 0.3cm

\nt \textbf{Main Theorem}: \emph{Every topological group is a
group retract of a minimal group.}

\vskip 0.3cm

It shows that Pestov's Conjecture is true. At the same time it
solves simultaneously Questions A and B
in a strong form. One of the conclusions is that the preservation
of minimality under quotients fails as strongly as possible when
passing from totally minimal to minimal groups. Note also that if
we do not require that $G$ is \emph{closed} in $M$ then this
weaker form of Question B (namely: every topological group $G$ is
a subgroup of a minimal group $M$) follows by a result of
Uspenskij \cite[Theorem 1.1]{Usp}. On the other hand in
Uspenskij's result the minimal group $M$ in addition is: a)
\emph{Raikov complete} (that is complete with respect to the two
sided uniformity); b) topologically simple and hence totally
minimal; c) \emph{Roelcke precompact} (that is the infimum
$\mathcal {U_L} \bigwedge \mathcal {U_R}$ (see Section
\ref{s:actions}) of right and left uniformities is precompact);
and d) preserves the weight of $G$.

 Our construction preserves some basic topological properties
like the weight, character, and the pseudocharacter. More
precisely: in the main theorem we prove that every topological
group $G$ can be represented as a group retract of a minimal group
$M$ such that simultaneously $w(M)=w(G)$, $\chi(M) = \chi(G)$ and
$\psi(M)=\psi(G)$ hold. In particular, if $G$ is metrizable (or
second countable) then the same is true for $M$. Moreover if $G$
is Raikov complete or \emph{Weil complete} (the latter means that
$G$ is complete with respect to the right uniformity) then in
addition we can assume that $M$ also has the same property. This
gives an immediate negative answer to the following

\ssk


\nt {\textbf{Question C:} (Arhangel'skii (see also \cite[Section
3.3D]{CHR})) \nl \emph{Let $M$ be a minimal group which is Raikov
complete. Must $\chi(G)=\psi(G)$ ? What if $M$ is Weil complete ?}

\ssk


Note that minimal topological groups with different $\chi(G)$ and
$\psi(G)$ (but without completeness assumptions) constructed
independently by
 Pestov \cite{pest-G-delta}, Shakhmatov \cite{Sh}
and Guran \cite{Gu} (see also \cite[Notes 7.7]{DPS}). This was one
of the motivations of Question C.

 From our main theorem we derive also that in fact every compact
homogeneous Hausdorf space admits a transitive continuous action
of a minimal group (see Corollary \ref{ex:arh-pr} below). This
means that minimality
makes no obstacle
in this setting. This fact negatively answers
the following

\ssk

\nt {\textbf{Question D:} (Arhangel'skii \cite[Problem VI.4]{Arh}
(see also \cite[Section 3.3G]{CHR})) \nl \emph{Suppose that a
minimal group acts continuously and transitively on a compact
Hausdorff space. Must $X$ be a dyadic space ? Must $X$ be a
Dugundji space ?}

\sk

In the proof of our main result we essentially use the methods of
\cite{me-min}. The main idea of \cite{me-min} was to introduce a
systematic method for constructing minimal groups using group
representations and generalized Heisenberg groups. For instance we
already proved (see \cite[Theorem 4.8]{me-min} or Theorems
\ref{t:BR} and \ref{t:BRapp} below) that if a group $G$ is
\emph{birepresantable}, that is if it admits sufficiently many
representations into continuous bilinear mappings (in short:
\emph{BR-group}; see Definition \ref{d:BR}.3), then $G$ is a group
retract of a minimal group. In the present paper we explore this
reduction by showing that in fact every topological group is a
BR-group. In the proof we use some new results about
representations into bilinear mappings. We show (see Theorem
\ref{t:matUC}) for instance that a bounded function $f: G \to \R$
on $G$ is left and right uniformly continuous if and only if $f$
is a \emph{matrix coefficient} $f=m_{v,\psi}$ (and hence
$f(g)=\psi(vg)$ for every $g \in G$) of a continuous Banach
co-representation $h: G \to Iso(V)$ by linear isometries such that
$v \in V$ and $\psi \in V^*$ is a $G$-continuous vector.
This result was inspired by a recent joint paper with Eli Glasner
\cite{GM2} (characterizing {\it strongly uniformly continuous}
functions on a topological group $G$ in terms of suitable matrix
coefficients). The technic in the latter result, as in some
related results of \cite{me-op, me-nz, GM2}, is based on a
dynamical modification of a celebrated factorization theorem in
Banach space theory discovered by Davis, Figiel, Johnson and
Pelczy\'nski \cite{DFJP}.

For the readers convenience, in the appendix (Section
\ref{s:appendix}) we include some proofs of \cite{me-min}.

\sk

\nt \textbf{Acknowledgment:} I would like to thank D. Dikranjan,
E. Glasner, M. Fabian, G. Itzkowitz, V. Pestov, V. Tarieladze and
V. Uspenskij for stimulating conversations and helpful
suggestions.


 \sk
\section{Preliminaries: actions and semidirect products}
\label{s:actions}
\sk


Let $X$ be a topological space. As usual denote by $w(X),
\chi(X),\psi(X), d(X)$ the weight, character, pseudocharacter and
the density of $X$ respectively.
All cardinals are assumed to be infinite.
%



 A (\emph{left}) \emph{action} of a topological group
$G$ on a space $X$, as usual, is a function $\pi : G \times X \to
X, \ \pi(g,x):=gx$ such that always $g_1(g_2x)=(g_1g_2)x$ and
$ex=x$ hold, where $e=e_G$ is the neutral element of $G$.
 Every $x \in X$ defines an \emph{orbit map} ${\tilde x}: G \to X, \ g
\mapsto gx$. Also every $g \in G$ induces a \emph{$g$-translation}
$\pi^g: X \to X, \ x \mapsto gx$. If the action $\pi$ is
continuous then we say that $X$ is a $G$-\emph{space}. Sometimes
we write it as a pair $(G,X)$.

Let $G$ act on $X_1$ and on $X_2$. A map $f: X_1 \to X_2$ is said
to be a \emph{$G$-map} if $f(gx)=gf(x)$ for every $(g,x) \in G
\times X_1$. A {\it $G$-compactification} of a $G$-space $X$ is a
continuous $G$-map $\alpha: X \to Y$ with a dense range
into a compact $G$-space Y.

 A {\it right action} $X \times G \to X$ can be defined analogously. If
$G^{op}$ is the {\it opposite group} of $G$ with the same topology
then the right $G$-space $(X,G)$ can be treated as a left
$G^{op}$-space $(G^{op},X)$ (and vice versa). A map $h: G_1 \to
G_2$ between two groups is a {\it co-homomorphism} (or, an
\emph{anti-homomorphism}) if $h(g_1g_2)=h(g_2)h(g_1)$. This
happens iff $h: G_1^{op} \to G_2$ (the same assignement) is a
homomorphism.

For a real normed space $V$ denote by $B_V$ its closed unit ball
$\{v \in V:  ||v|| \leq 1 \}$.
Denote by $Iso(V)$ the topological group of all linear (onto)
isometries $V \to V$ endowed with the \emph{strong operator
topology}. This is the topology of pointwise convergence inherited
from $V^V$. Let $V^*$ be the dual Banach space of $V$ and let
$$
<,>: V \times V^* \to \R, \ \ (v,\psi) \mapsto <v,\psi>=\psi(v)
$$
be the canonical (always continuous) bilinear mapping. Let $\pi: G
\times V \to V$ be a continuous left action of $G$ on $V$ by
linear isometries. This is equivalent to saying that the natural
homomorphism $h: G \to Iso(V), \ g \mapsto \pi^g$ is continuous.
The \emph{adjoint action} $G \times V^* \to V^*$ is defined by
$g\psi(v):=\psi(g^{-1}v)$. Then the corresponding canonical form
is $G$-invariant. That is
$$
<gv,g\psi>=<v,\psi> \ \ \  \forall \ (g,v,\psi) \in G \times V
\times V^*.
$$

Similarly, if $V \times G \to V$ is a continuous \emph{right}
action of $G$ on $V$ by linear isometries. Then the corresponding
\emph{adjoint action} (from the left) $G \times V^* \to V^*$ is
defined by $g\psi(v):=\psi(vg)$. Then we have the following
equality
$$
<vg,\psi>=<v,g\psi> \ \ \  \forall \ (g,v,\psi) \in G \times V
\times V^*.
$$

Adjoint actions of $G$ on $V^*$ does not remain continuous in
general (see for example \cite{me-fr}).

%
%

The Banach algebra (under the supremum norm) of all continuous
real valued bounded functions on a topological space $X$ will be
denoted by $C(X)$.
 Let $(G,X)$ be a left $G$-space. It induces the right
action $C(X)\times G \to C(X)$, where $(fg)(x)=f(gx)$, and the
corresponding co-homomorphism $h: G \to Iso(C(X))$.
%
While the $g$-translations $C(X) \to C(X)$ are continuous (being
isometric), the orbit maps ${\tilde f}: G \to C(X), \ g \mapsto
fg$ are not necessarily continuous. However if $X$ is a compact
$G$-space then every ${\tilde f}$ is continuous and equivalently
the action $C(X) \times G \to C(X)$ is continuous.

For every topological group $G$ denote by $RUC(G)$ the Banach
subalgebra of $C(G)$ of {\em right uniformly continuous\/} (some
authors call these functions {\it left uniformly continuous})
bounded real valued  functions on $G$. These are the functions
which are uniformly continuous with respect to the {\it right
uniform structure} $\mathcal {U_R}$ on $G$. Thus, $f\in RUC(G)$
iff for every $\eps>0$ there exists a neighborhood $V$ of the
identity element $e\in G$ such that $\sup_{g\in
G}|f(vg)-f(g)|<\eps$ for every $v \in V$. It is equivalent to say
that the orbit map $G \to C(G), \hskip 0.2cm g \mapsto fg$ is norm
continuous where $fg$ is the {\it left translation} of $f$ defined
by $(fg)(x):=f(gx)$.

Analogously can be defined the algebra $LUC(G)$ of {\it left
uniformly continuous} functions (and the {\it right
translations}). These are the functions which are uniformly
continuous with respect to the {\it left uniform structure}
$\mathcal {U_L}$ on $G$.

Denote by $\mathcal {U_L} \bigwedge \mathcal {U_R}$ the
\emph{lower uniformity} of $G$. It is the infimum (greatest lower
bound) of left and right uniformities on the set $G$. The
intersection $UC(G):=RUC(G) \cap LUC(G)$ is a left and right
$G$-invariant closed subalgebra of $RUC(G)$. Clearly, for every
bounded function $f: G \to \R$ we have $f \in UC(G)$ iff $f:
(G,\mathcal {U_L} \bigwedge \mathcal {U_R}) \to \R$ is uniformly
continuous.
We need the following important and non-trivial fact.

\begin{lem} \label{l:Roelcke}
\ben
\item \emph{(Roelcke-Dierolf \cite{RD})}
For every topological group $G$ the lower uniformity $\mathcal
{U_L} \bigwedge \mathcal {U_R}$ generates the given topology of
$G$.
\item
For every topological group $G$ the algebra $UC(G)$
separates points from closed subsets in $G$. \een
\end{lem}
\begin{proof} (1) \
See Roelcke-Dierolf \cite[Proposition 2.5]{RD}.

(2) Follows from (1).
\end{proof}

Note that in general the infimum $\mu_1 \bigwedge \mu_2$ of two
compatible uniform structures on a topological space $X$ need not
be compatible with the topology of $X$
(see for example \cite[Remark 2.2]{Tk}).

\sk

Let $(X,\t)$ and $(G,\s)$ be topological groups and
$$
\a: G \times X \to X, \ \ \a(g,x)=gx=g(x)
$$
 be a given (left) action. We say that $X$
is a \emph{$G$-group} if $\a$ is continuous and every
$g$-translation $\a^g: X \to X$ is a group automorphism of $X$.
For every $G$-group $X$ denote by $X \leftthreetimes_{\a} G$ the
corresponding \emph{topological semidirect product} (see for
example \cite[Section 6]{RD} or \cite[Ch. 7]{DPS}). As a
topological space this is the product $X \times G$. The standard
multiplicative group operation is defined by the rule: for a pair
$(x_1,g_1)$, $(x_2,g_2)$ in $X \leftthreetimes_{\a} G$ let
$$
(x_1,g_1) \cdot (x_2,g_2):= (x_1\cdot g_1(x_2), \ g_1 \cdot g_2).
$$
 Sometimes the closed normal subgroup $X \times \{e_G\}$ of $X
\leftthreetimes_{\a} G$ will be identified with $X$ and similarly,
the closed subgroup $\{e_X\} \times G$ will be identified with
$G$.
The projection $p: X \leftthreetimes_{\a} G \to G, \ p(x,g)=g$ is
a group homomorphism and also a retraction.
In particular, $G$ is a quotient of $X \leftthreetimes_{\a} G$.
The kernel of this retraction $ker(p)$ is just $X \times \{e_G\}$.

\begin{defn} \label{d:t-exactness} 
Let $(X,\t)$ and $(G,\s)$ be Hausdorff topological groups and
$$
\a: G \times X \to X, \ \ \a(g,x)=gx=g(x)
$$
 be a continuous action by group automorphisms.
\ben
\item
 The action $\a$ is {\it topologically exact} (\emph{t-exact}, for short)
if there is no strictly coarser, \emph{not necessarily Hausdorff},
group topology $\s' \varsubsetneq \s$ on $G$ such that $\a$ is
$(\s',\t,\t)$-continuous (See \cite{me-Gmin}).
\item
More generally, let $ \{\a_i: G \times Y_i \to Y_i\}_{i \in I}\ $
be a system of continuous $G$-actions on the groups $Y_i$. We say
that this system is \emph{t-exact} if there is no strictly
coarser, not necessarily Hausdorff, group topology $\s'
\varsubsetneq \s$ on $G$ such that all given actions remain
continuous.
\item $X$ is a \emph{$G$-minimal group} if there is no strictly coarser Hausdorff
group topology $\t' \varsubsetneq \t$ on $X$ such that $\a$
remains continuous with respect to the triple $(\s, \t', \t')$ of
topologies (See \cite{RS}).
\een
\end{defn}


\begin{lem} \label{l:G-min}
Let $(X \leftthreetimes_{\a} G, \g)$ be a topological semidirect
product. Suppose that $X$ is $G$-minimal with respect to $\a$.
Then for every coarser Hausdorff group topology $\g_1 \subset \g$
we have $\g_1|_X=\g|_X$.
\end{lem}
\begin{proof}
Since $P:=(X \leftthreetimes_{\a} G, \g_1)$ is a topological group
the conjugation map $$(P, \g_1) \times (P, \g_1) \to (P, \g_1), \
(a,b) \to aba^{-1}$$ is continuous. Then its restriction $$(G,
\g_1|_G) \times (X, \g_1|_X) \to (X, \g_1|_X), \ (g,x) \to
g(x)=gxg^{-1}$$ is also continuous. Since $\g_1|_G \subset \g|_G$
it follows that the action of the given group $(G,\g|_G)$ on the
Hausdorff group $(X, \g_1|_X)$ is continuous, too. Since $\g_1|_X
\subset \g|_X$ and $X$ is $G$-minimal we obtain $\g_1|_X=\g|_X$.
\end{proof}

\begin{thm} \label{t:semi-products}
Let $G$ be a Hausdorff topological group and let $X$ be an abelian
$G$-minimal group. Then if the given action $\a: G \times X \to X$
is t-exact then $X \leftthreetimes_{\a} G$ is minimal.
\end{thm}
\begin{proof}
See Theorem 1.4 of \cite{me-Gmin}
(or combine Lemma \ref{l:G-min} and Corollary \ref{c:2.8}).
\end{proof}

\begin{lem} \label{l:completionofsemid}
\cite[Proposition 12.5]{RD} Let $P:=X \leftthreetimes_{\a} G$ be a
topological semidirect product. If $X$ and $G$ are Raikov complete
(Weil complete) then $P$ is a Raikov complete (resp., Weil
complete) group.
\end{lem}

\begin{remarks} \label{r:definitionsand}
\ben
\item
Note that in \cite{me-min} the original definition of t-exactness
contains a superfluous condition of \emph{algebraic exactness.}
The latter means that the kernel of the action $ker(\a):=\{g \in
G: \ gx=x \ \ \forall \ x \in X \}$ is trivial. The reason is that
since $G$ is Hausdorff every t-exact action is algebraically
exact. Indeed assuming the contrary let $H:=ker(\a)$ be the
\emph{nontrivial} kernel of the action $\a$. Consider the quotient
group $G/H$ with the coset topology $\tau/H$ and the map $q: G \to
G/H, \ g \mapsto q(g)=gH$. Then the induced action of $G/H$ on $X$
is continuous. It follows that the preimage topology
$\t':=q^{-1}(\tau/H)$ on $G$ is a group topology such that $\a$
remains continuous. Since $q^{-1}(\tau/H)$, being not Hausdorff,
is strictly coarser than the original (Hausdorff) topology of $G$
we obtain that $\a$ is not t-exact.
\item
Suppose that $X$ is a $G$-group under the action $\pi: G \times X
\to X$ such that the semidirect product $X \leftthreetimes_{\pi}
G$ is minimal. Then if $\pi$ is algebraically exact then $\pi$
necessarily is t-exact. Indeed, otherwise there exists a strictly
coarser group topology $\t'$ on $G$ such that $\a: (G,\t') \times
X \to X$ remains continuous. Since $X$ is Hausdorff for every $x
\in X$ and every $g$ from the $\t'$-closure $cl_{\t'}(\{e\})$ of
the singleton $\{e\}$ we have $gx=x$. Since the action is
algebraically exact we get $cl_{\t'}(\{e\})=\{e\}$.
Thus, $\t'$ is Hausdorff. Then the semidirect product $X
\leftthreetimes_{\a} (G,\t')$ is a Hausdorff topological group and
its topology is strictly coarser than the original topology on $X
\leftthreetimes_{\a} (G, \t)$. This contradicts the minimality of
the latter group.
\item
The direct product $X \times G$ of two minimal abelian (even
cyclic) groups $X$ and $G$ may not be minimal. Take for example
$X=G=(\Z, \t_p)$ with the $p$-adic topology $\t_p$ (see Do\"\i
chinov \cite{Do}). Since $X$ is minimal it also can be treated as
a $G$-minimal group \wrt the trivial action of $G$ on $X$. Then
the direct product is just the semidirect product in our setting.
It follows (as expected, of course) that the t-exactness is very
essential in Theorem \ref{t:semi-products}. This example also
demonstrates that the quote "not necessarily Hausdorff" in the
Definition \ref{d:t-exactness}.2 of the t-exactness cannot be
omitted.
\item
If $X$ is a locally compact Hausdorff group and $G$ is a subgroup
of $Aut(X)$ endowed with the standard {\it Birkhoff topology} (see
\cite{DPS,me-min}) then the corresponding action is t-exact.
\item
For every normed space $V$ and a topological subgroup $G$ of
$Iso(V)$ the action of $G$ on $V$ is t-exact.
\item
According to \cite[Example 10]{EDS} there exists a totally minimal
precompact group $X$ such that a certain semidirect product $X
\leftthreetimes \Z_2$ with the two-element cyclic group $\Z_2$ is
not minimal. The given action of $\Z_2$ on $X$ is t-exact. Indeed,
by the construction the action is not trivial. On the other hand
every strictly coarser group topology on the (discrete) group
$\Z_2$ is the trivial topology. This example demonstrates that
Theorem \ref{t:semi-products} is not true in general for
non-abelian $X$.
 \een
\end{remarks}


\sk
\section{Generalized Heisenberg groups}
\label{s:GH}
\sk

 We need a natural generalization of the classical
three dimensional Heisenberg group. This generalization is based
on semidirect products defined by biadditive mappings. See, for
example, \cite{Re, Mi, me-min}). For additional new properties and
applications of this construction we refer also to \cite{me-Gmin,
me-hei, DM}.

Let $E,F,A$ be abelian groups. A map $w: E \times F \to A$ is said
to be {\it biadditive} if the induced mappings
$$
\om_x: F \to A, \ w_f: E \to A, \ \ \om_x(f):=\om(x,f)=:\om_f(x)
$$
are homomorphisms for all $x \in E$ and $f \in F$. Sometimes we
look at the elements f of F as functions defined on E, i.e., the
value f(x), for an element x of E is defined as $\omega(x,f)$. We
say that $\om$ is {\it separated} if the induced homomorphisms
separate points. That is, for every $x_0 \in E, f_0 \in F$ there
exist $f \in F, x \in E$ such that $f(x_0) \neq 0_A, f_0(x) \neq
0_A$, where $0_A$ is the zero element of $A$.

\begin{defn}
Let $E, F$ and $A$ be
abelian topological groups and $\om: E \times F \to A$ be a
continuous biadditive mapping. Denote \footnote{The notation of
the present paper about generalized Heisenberg groups and some
related objects not always agree with the notation of
\cite{me-min}. } by
$$
H(\om)= (A \times E) \leftthreetimes_{\om^{\nabla}} F
$$
the semidirect product (say, {\it generalized Heisenberg group}
induced by $\om$) of $F$ and the group $A \times E$ with respect
to the action
$$
\om^{\nabla}: F \times A \times E \to A \times E, \ \  \
f(a,x)=(a+f(x),x).
$$
The resulting group, as a topological space, is the product $A
\times E \times F $.
The group operation is defined by the following rule:
for a pair
$$u_1=(a_1,x_1,f_1), \hskip 0.4cm u_2=(a_2,x_2,f_2)$$
define
$$u_1 \cdot u_2 = (a_1+a_2+f_1(x_2), x_1+x_2,
f_1 +f_2).$$ 
\end{defn}

\sk

Then $H(\om)$ becomes a two-step nilpotent (in particular,
metabelian)
topological group.

%
%
Elementary computations for the commutator $[u_1,u_2]$ give
$$[u_1,u_2] = u_1u_2u_1^{-1}u_2^{-1}= (f_1(x_2)-f_2(x_1),0_E,0_F).$$

If $\om$ is separated then the center of the group $H(\om)$ is the
subgroup $A$.

%


The very particular case of the canonical bilinear form $<,>: \R^n
\times \R^n \to \R$ defines the classical 2n+1-dimensional
Heisenberg group. If $A$, $E$ and $F$ are complete groups then by
Lemma \ref{l:completionofsemid} the corresponding Heisenberg group
$H(\om)$ is Weil complete.

In \cite{me-min} we show that generalized Heisenberg groups are
useful in the theory of minimal topological groups. One of the
results obtained there is: for every \lca group $G$ the Heisenberg
group $H(\Delta)= \T \times G^* \leftthreetimes_{\Delta} G$ of the
canonical biadditive mapping $\Delta: G^* \times G \to \T, \hskip
0.2cm \Delta(\chi,g)= \chi(g)$, where $\T$ denotes the circle
group, is minimal. It follows that every \lca group is a group
retract of a locally compact minimal group.

For more examples of minimal groups that come from biadditive
mappings see \cite{DM}.

The following definition generalizes slightly \cite[Definition
3.1]{me-min}.

\sk

\begin{defn} \label{d:str-duality}

 Let $E$ and $F$ be (semi)normed spaces
and $\om: E \times F \to \R$ be a bilinear map. We say that:
 \ben
\item
 $\om$ is a \emph{left strong duality} if for every norm-unbounded
sequence $x_n \in E$ the subset $\{f(x_n): n \in \N, \ f \in B_F
\}$ is unbounded in $\R$. This is equivalent to saying that
$\{f(x_n): \ n\in \N, \ f \in O \} = \R$ for every \nbd $O$ of the
zero $0_F$ (or, even, for every nonempty open subset $O$ in $F$).
\item $\om$ is a a \emph{right strong duality} if for every norm-unbounded
sequence $f_n \in F$ the subset $\{f_n(x): \ n \in \N, \ x \in B_E
\}$ is unbounded in $\R$.
\item
left and right strong duality we call simply a \emph{strong
duality}.
 \een
\end{defn}

If $\om: E \times F \to \R$ is a strong duality with normed spaces
$E$ and $F$ then $\om$ necessarily is separated. Indeed, let $v
\in V, \ v \neq 0_E$ and $f(v)=0$ for every $f \in F$. Since $E$
is a normed space and $v \neq 0_E$ we have $\norm{v}
>0$. Then the sequence $x_n:=nv_n$ is unbounded in $V$. On the
other hand $\{f(x_n): \ n\in \N, f \in B_F \}=\{0\}$. This means
that $\om$ is not left strong. Similarly can be proved the case of
$f \in F, \ f \neq 0_F$ with $f(v)=0$ for every $v \in V$.

\begin{example} \cite{me-min}
\ben
\item
 For every normed space $V$ the canonical bilinear form $<,>: V
\times V^* \to \R$ is a strong duality.
\item
For every locally compact group $G$ the natural bilinear form
$$\om: L^1(G) \times \mathcal{K}(G) \to \R$$ is a strong duality.

Here $\mathcal{K}(G)$ is the normed space of all continuous real
valued functions with compact supports endowed with the sup-norm.
It can be treated as a \emph{proper} subspace of
$L^1(G)^*:=L^{\infty}(G)$ such that $\om$ is a restriction of the
canonical form $L^1(G) \times L^1(G)^* \to \R$ (hence the second
example is not a particular case of (1)). \een
\end{example}

Let $\om: E \times F \to \R$ be a
separated
bilinear mapping. Then the generalized Heisenberg group $H(\om)=\R
\times E \leftthreetimes_{\om^{\nabla}} F^*$ is not minimal.
Indeed the center of a minimal group must be minimal (see for
example \cite[Proposition 7.2.5]{DPS}) and the center of $H(\om)$
is the subgroup $\R$ which is not minimal. Note however that the
subgroups $V$ and $V^*$ are \emph{relatively minimal} in $H(\om)$
(see \cite{DM, me-hei}) for the canonical duality $\om=<,>: V
\times V^* \to \R$ for every normed space $V$.

Now we define as in \cite{me-min} the semidirect product
$$H_+(\om):=H(\om) \leftthreetimes_{\a} R_+ = (\R \times E
\leftthreetimes_{\om^{\nabla}} F^*) \leftthreetimes_{\a} \R_+$$
where $\R_+$ is the multiplicative group of all positive reals and
$\a$ is the natural action
$$\a: \R_+ \times H(\om) \to H(\om), \ \
\a(t,(a,c,f))=(ta,tx,f).$$ Observe that the third coordinate after
the $t$-translation is just "$f$" and not "$tf$".

It turns out that $H_+(\om)$ is minimal under natural restrictions
providing a lot of examples of minimal groups.

\begin{thm} \label{str-dual-ismin} \cite{me-min}
For every strong duality $\om: E \times F \to \R$ with normed
spaces $E$ and $F$ the corresponding group $H_+(\om)=H(\om)
\leftthreetimes_{\a} \R_+$ is minimal.
\end{thm}
\begin{proof}
See \cite[Proposition 3.6]{me-min} or Theorem \ref{t:plus} below.
\end{proof}

Furthermore, using (twice) Lemma \ref{l:completionofsemid} we get

\begin{lem} \label{l:complconcrete}
For every continuous bilinear mapping $\om: E \times F \to \R$
with Banach spaces $E$ and $F$ the corresponding group
$H_+(\om)=H(\om) \leftthreetimes_{\a} R_+=(\R \times E
\leftthreetimes_{\om^{\nabla}} F^*) \leftthreetimes_{\a} \R_+$ is
Weil complete.
\end{lem}

%
%
%
%
%
%
%
%
%
%
%
%

\sk

\sk
\section{Group representations in bilinear forms}
\sk

Let $G$ be a topological group. A
\emph{representation} (\emph{co-representation}) of $G$ on a
normed space $V$, is a
homomorphism (resp. co-homomorphism) $h: G \to Iso(V)$. Sometimes
we give the representation (co-representation) by the
corresponding linear isometric left (resp. right) action $G \times
V \to V, \ (g,v) \mapsto gv=h(g)(v)$ (resp., $V \times G \to V, \
(v,g) \mapsto vg=h(g)(v)$).

\begin{defn} \label{d:mat}
  Let $V, W$ be normed spaces and $\om: V \times W \to \R$
be a continuous bilinear map.
\ben
\item Let $h_1: G \to Iso(V)$ and $h_2: G \to
Iso(W)$ both are homomorphisms. We say that the pair $(h_1,h_2)$
is a \emph{birepresentation} of $G$ in $\om$ if $\om$ is
$G$-invariant. That is, $$\om(gv,g\psi)=\om(v,\psi) \ \ \quad
\quad \forall
 \ (g,v,\psi) \in G \times V \times W.$$
\item
Let $h_1: G \to Iso(V)$ be a co-homomorphism and and $h_2: G \to
Iso(W)$ be a homomorphism. We say that the pair $(h_1,h_2)$ is a
\emph{co-birepresentation }of $G$ in $\om$ if
$$\om(vg,\psi)=\om(v,g\psi) \ \ \quad \quad \forall
 \ (g,v,\psi) \in G \times V \times W.$$
 \een
\end{defn}


\sk

The following definition is one of the key ideas of \cite{me-min},
as well as of the present paper. Assume that the pair $\a_1: G
\times E \to E, \ \a_2: G \times F \to F$ is a birepresentation
$\Psi$ of $G$ in $\om: E \times F \to \R$. By the \emph{induced
group} $M_+(\Psi)$ of the given birepresentation $\Psi$, we mean
the topological semidirect product $H_+(\om) \leftthreetimes_{\pi}
G$, where the action
$$\pi: G
\times H_+(\om) \to H_+(\om)$$ is defined by $\pi(g,
(a,x,f,t)):=(a,gx,gf,t)$, where $gx=\a_1(g,x)$ and $gf=\a_2(g,f)$.

More generally, let
$$\Phi:=\{\Phi_i\}_{i \in I}=\{\om_i: E_i
\times F_i \to \R, \  \ \a_{1i}: G \times E_i \to E_i, \ \
\a_{2i}: G \times F_i \to F_i\}_{i \in I}$$ be a system of
continuous $G$-birepresentations. By the \emph{induced group }
$M_+(\Phi)$ of the system $\Phi$ we mean the semidirect product
$$
\prod_{i \in I} M_+(\om_i) \leftthreetimes_{\pi} G
$$
where the action $$\pi: G \times \prod_{i \in I} M_+(\om_i) \to
\prod_{i \in I} M_+(\om_i)$$ is defined coordinatwise by means of
the following system $\{\pi_i\}_{i \in I}$ of actions (defined
above):

$$\pi_i: G
\times H_+(\om_i) \to H_+(\om_i), \quad \pi_i(g,
(a,x,f,t)):=(a,gx,gf,t).$$

\begin{remark} \label{r:opp} Let $h: G \to Iso(V)$ be a
 co-homomorphism (homomorphism). Denote by $h^{op}$ the \emph{associated
homomorphism} (resp. co-homomorphism) $h^{op}: G \to Iso(V), \ g
\mapsto h(g^{-1}).$ Then the pair $(h_1,h_2)$ is a
co-birepresentation of $G$ iff $(h_1^{op},h_2)$ is a
birepresentation of $G$ in $\om$.

Analogously for every system $\Phi$ of co-birepresentations we can
define the naturally associated system $\Phi^{op}$ of
birepresentations.
\end{remark}

\begin{defn} \label{d:BR}
 (1) \ let $\Phi:=\{\om_i: E_i \times F_i \to \R, \ \
\a_{1i}: G \times E_i \to E_i, \ \a_{2i}: G \times F_i \to
F_i\}_{i \in I}$ be a system of continuous $G$-birepresentations.
We say that this system is \emph{topologically exact}
(\emph{t-exact}) if $(G,\tau)$ is a Hausdorff topological group
and for every strictly coarser (not necessarily Hausdorff) group
topology $\tau' \subset \tau$ on $G$ there exists an index $i \in
I$ such that one of the actions $\a_{1i}: (G,\tau') \times E_i \to
E_i$ or $\a_{1i}: (G,\tau') \times F_i \to F_i$ is not continuous.

Similarly can be defined t-exact systems of
\emph{co-birepresentations}. Obviously the system of
co-birepresentations is t-exact iff the associated system (see
Remark \ref{r:opp}) of birepresentations is t-exact.


(2) \ \cite[Definition 4.7]{me-min} We say that a topological
group $G$ is \emph{birepresentable} (shortly: \emph{BR-group}) if
there exists a t-exact system $\Phi$ of linear birepresentations
(equivalently, co-birepresentations (see Remark \ref{r:opp})) of
$G$ in \emph{strong dualities} $\om_i$.
\end{defn}

The following theorem from \cite{me-min} (See also Theorem
\ref{t:BRapp} below) is one of the crucial results in our setting.

\begin{thm} \label{t:BR} \cite[Theorem 4.8]{me-min}
\ben
\item
 Let $\Phi$ be a t-exact system of
 $G$-birepresentations into strong dualities
 $\om_i: E_i \times F_i \to \R$ with normed spaces $E_i$ and $F_i$.
 Then the corresponding induced group
$M_+(\Phi)$ is minimal.
\item
For every BR-group $G$ there exists a continuous group retraction
$p: M \to G$ such that $M$ and also the kernel $ker(p)$ are
minimal.
 \een
\end{thm}

Some results of \cite{me-min, me-fr} show that many important
groups (like additive subgroups of locally convex spaces and
locally compact groups) are BR-groups. One of the main results of
the present paper (see Theorem \ref{t:allBR}) shows that in fact
every topological group is a BR-group. Furthermore, in the
Definition \ref{d:BR}.2 we can always choose a system with
$|I|=1$; that is a system $\Phi$ with a single birepresentation.

 \sk
\section{Matrix Coefficients of group representations}
\sk

We generalize the usual notion of matrix coefficients to the case
of arbitrary (not necessarily canonical) bilinear mappings.

\begin{defn} \label{d:mat}
  Let $V, W$ be (semi)normed spaces and let $\om: V \times W \to \R$
be a continuous bilinear map. Let $h_1: G \to Iso(V)$ be a
co-representation and $h_2: G \to Iso(W)$ be a representation such
that the pair $(h_1,h_2)$ is a co-birepresentation of $G$ in
$\om$. \ben
\item
For every pair of vectors $v \in V$ and $\psi \in W$ define the
\emph{matrix coefficient} $ m_{v, \psi}$ as the following function
$$
m_{v, \psi}: G \to \R, \hskip 0.3cm g \mapsto \psi(vg)$$ (where
$\psi(vg)=\om(vg,\psi)=<vg, \psi>=<v, g \psi>$).
\item
We say that a vector $v \in V$ is \emph{$G$-continuous} if the
corresponding orbit map $\tilde{v}: G \to V, \ \tilde{v}(g)=vg$,
defined through $h_1: G \to Iso(V)$, is norm continuous. Similarly
one can define a $G$-\emph{continuous} vector $\psi \in W$.
\item
We say that a matrix coefficient $m_{v,\psi}: G \to \R$ is
\emph{bicontinuous} if $v \in V$ and $\psi \in W$ are
$G$-continuous vectors. If $(h_1,h_2)$ is a continuous
co-birepresentation then every corresponding matrix coefficient is
bicontinuous.
 \een
\end{defn}

First we need the following

\begin{lem} \label{l:cont}
Let $ \omega: V \times W \to \R, \ (v,f) \mapsto
\om(v,f)=<v,f>=f(v) $ be a bilinear mapping defined for
(semi)normed spaces $V$ and $W$. \TFAE
 \ben
\item
$\om$ is continuous.
\item
For some constant $c >0$ the inequality $|<v,f>| \leq c \cdot
\norm{v} \cdot \norm{f}$ holds for every $(v,f) \in V \times W$.
 \een
\end{lem}
\begin{proof} (1) $\Longrightarrow$ (2): \
By the continuity of $\om$ at $(0_V, 0_W)$ there exists a constant
$\eps
>0$ such that the inequalities $\norm{v} \leq \eps, \ \norm{f} \leq \eps$ imply
$|f(v)| \leq 1$. Then $|(\eps \cdot \frac{f}{\norm{f}})(\eps \cdot
\frac{v}{\norm{v}})| \leq 1$ holds for every $(v,f) \in V \times
W$. It follows that $|f(v)| \leq \frac{1}{\eps^2} \cdot \norm{v}
\cdot \norm{f}$ for every $(v,f) \in V \times W$.

(2) $\Longrightarrow$ (1): \ Is trivial.
\end{proof}

\begin{defn} \label{d:regular}
We say that a bilinear mapping $\om: V \times W \to \R$ is
\emph{regular} if
$$|<v,f>| \leq \norm{v} \cdot \norm{f}$$ holds for every $(v,f)
\in V \times W$. If $\om$ in addition is a strong duality then we
call it a \emph{regular strong duality}.
\end{defn}

For example, the canonical bilinear mapping $<,>: V \times V^* \to
\R$ for every normed space $V$ (and hence its any restriction) is
regular. Every regular bilinear mapping is continuous by Lemma
\ref{l:cont}.

\sk

The following observation is a modification of \cite[Fact
3.5.2]{me-nz}.

\begin{lem} \label{introversion}
Let $\om: V \times W \to \R$ be a continuous bilinear mapping.
Assume that the pair $h_1: G \to Iso(V), \ h_2: G \to Iso(W)$ is a
(not necessarily continuous) co-birepresentation of $G$ in $\om$.
Then every bicontinuous matrix coefficient $f=m_{v,\psi}: G \to
\R$ is left and right uniformly continuous on $G$ (that is, $f \in
UC(G)$).
\end{lem}
\begin{proof}
Since $\om$ is continuous by Lemma \ref{l:cont} there exists a
constant $c >0$ such that
$$
|<x,y>| \leq c \cdot \norm{x} \cdot \norm{y} \quad \forall \ (x,y)
\in V \times W.
$$
Since $h_1(G) \subset Iso(V)$ and $h_2(G) \subset Iso(W)$, we have
$\norm{xg} = \norm{x}$ and $\norm{gy} = \norm{y}$ for every
$(g,x,y) \in G \times V \times W$. Clearly
$|m_{v,\psi}(g)|=<vg,\psi> \leq c \cdot \norm{v} \cdot
\norm{\psi}$. Hence $m_{v,\psi}$ is a bounded function. In order
to establish that $f=m_{v,\psi} \in LUC(G)$, observe that
$$|f(gu)-f(g)|=|m_{v,\psi}(gu) - m_{v,\psi}(g)|=| <vgu, \psi> - <vg,
\psi>|=$$
$$ |<vg,u \psi> - <vg, \psi>|
\leq c \cdot ||vg|| \cdot ||u \psi-\psi|| = c \cdot \norm{v} \cdot
||u\psi-\psi||.$$ Here, using the $G$-continuity of the vector
$\psi$ in $W$, we get that $f \in LUC(G)$.

Similar verification is valid for the second case $f \in RUC(G)$.
\end{proof}

Principal Theorem \ref{t:matUC} below shows that the
representability of a function as a bicontinuous matrix
coefficient in fact characterizes functions from $UC(G)=LUC(G)
\cap RUC(G)$.

\sk

\begin{defn} \label{d:local-separ}
We say that a family $S$ of continuous functions on a topological
group $G$ is a \emph{local separating family} if $S$ separates the
identity $e \in G$ from the closed subsets of $G$ that do not
contain the identity. That is, for every neighborhood $U$ of $e$
in $G$ there exist: $f \in S$, $\eps >0$ and a real number $r \in
\R$ such that $f(e)=r$ and $f^{-1}(r-\eps, r+\eps) \subset U$.
\end{defn}

\begin{lem} \label{l:separation}
Let $(G,\tau)$ be a topological group and $S$ be a local
separating family of functions. Suppose that $\tau' \subset \tau$
is a coarser group topology on $G$ such that every $f \in S$ is
continuous on a topological group $(G, \tau')$. Then $\tau'
=\tau$.
\end{lem}
\begin{proof} Observe that by our assumption the homomorphism of groups $1_G:
(G,\tau') \to (G,\tau), \ g \mapsto g$ is continuous at the
identity $e$. Hence this homomorphism is continuous. This implies
that, $\tau \subset \tau'$. Hence, $\tau' =\tau$, as required.
\end{proof}

\begin{lem} \label{l:Mat}
Let $(G,\tau)$ be a Hausdorff topological group and let
$$\Phi:=\{\om_i: E_i \times F_i \to \R,  \ \ \a_{1i}: E_i \times G
\to E_i, \ \a_{2i}: G \times F_i \to F_i\}_{i \in I}$$ be a system
of continuous co-birepresentations of $G$ into bilinear mappings
$\om_i$. Let
$$
{\mathcal M}_{\Phi}:=\{m_{v,\psi}: G \to \R:  \ \ \ \ (v,\psi) \in
E_i \times F_i\}_{i \in I}
$$
be the family of corresponding matrix coefficients. Suppose that
${\mathcal M}_{\Phi}$ is a local separating family on $G$. Then
$\Phi$ is t-exact.
\end{lem}
\begin{proof} 
%
Assume that $\tau_1 \subset \tau$ is a coarser group topology on
$G$ such that all given co-birepresentations are still continuous.
Then by Lemma \ref{introversion}, every matrix coefficient
$m_{v,\psi}: (G,\tau_1) \to \R$ is (uniformly) continuous for each
$m_{v,\psi} \in {\mathcal M}_{\Phi}$. By our assumption ${\mathcal
M}_{\Phi}$ is a local separating family. By Lemma
\ref{l:separation} we get $\tau_1=\tau$. This means that the
system $\Phi$ is t-exact.
\end{proof}

\sk


The following definition was inspired by \cite{GM2}. Namely by the
concept of \emph{Strong Uniform Continuity} (SUC).

\begin{defn} \label{d:small}
%
Let $h: G \to Iso(V)$ be a continuous co-representation on a
normed space $V$ and $x_0 \in V^*$. We say that a
subset $M \subset V$ is \emph{SUC-small at} $x_0$ if
for every $\eps
>0$ there exists a \nbd $U$ of $e$ such that
$$sup_{v \in M} |<v, ux_0> -<v, x_0>| \leq \eps \ \ \quad \forall \ u \in
U.$$

\end{defn}

%

\ssk

We collect here some useful properties of SUC-smallness.

\begin{lem} \label{l:small}
Let $h: G \to Iso(V)$ be a continuous co-representation and $x_0
\in V^*$.
 \bit
\item [(a)] The family of SUC-small sets at $x_0$ is closed under
taking: subsets, norm closures, finite linear combinations and
convex hulls.
\item [(b)]
If $M_n \subset V$ is SUC-small at $x_0 \in V^*$ for every $n \in
\N$ then so is $\bigcap_{n \in \N} (M_n + \delta_n B_V)$ for every
positive sequence $\delta_n$ such that $\lim \delta_n =0$.
\item [(c)] Let $h: G \to Iso(V)$ be a
co-representation. For every $\psi \in V^*$ the following are
equivalent:
 \bit
\item [(i)]
The orbit map $\widetilde{\psi}: G \to V^*$ is norm continuous.
\item [(ii)]
$\mathcal B$ is SUC-small at $\psi$, where $\mathcal
B:=\{\breve{v}: V^* \to \R, \ x \mapsto \check{v}(x):=<v,x> \}_{v
\in B_V}$. \eit
\item [(d)] Let $h_2: G \to Iso(E)$ be a continuous co-representation
and let $\gamma: V \to E$ be a linear continuous $G$-map
(of right $G$-spaces). Assume that $M \subset E$ is an SUC-small
set at $\psi \in E^*$. Then $\gamma^{-1}(M) \subset V$ is
SUC-small at $\gamma^*(\psi) \in V^*$, where $\g^*: E^* \to V^*$
is the adjoint of $\g$. \eit
\end{lem}
\begin{proof}
Assertion (a) is straightforward.

(b): \ We have to show that the set $\bigcap_{n \in \N} (M_n +
\delta_n B_V)$ is SUC-small at $x_0$. Let $\eps
>0$ be fixed. Since $Gx_0$ is a bounded subset of $V^*$ one can
choose $n_0 \in \N$ such that $|v(gx_0)| < \frac{\eps}{4}$ for
every $g \in G$ and every $v \in \delta_{n_0} B_V$. Since
$M_{n_0}$ is SUC-small at $x_0$ we can choose a \nbd $U(e)$ such
that $|m(ux_0)-m(x_0)| < \frac{\eps}{2}$ for every $u \in U$ and
every $m \in M_{n_0}$. Now every element $w \in \bigcap_{n \in \N}
(M_n + \delta_n B_V)$ has a form $w=m + v$ for some $m \in
M_{n_0}$ and $v \in \delta_{n_0} B_V$. Then for every $u \in U$ we
have
$$|w(ux_0) - w(x_0)| \leq |m(ux_0)- m(x_0)| + |v(ux_0)| + |v(x_0)|
< \frac{\eps}{2} + \frac{\eps}{4}+ \frac{\eps}{4} = \eps.
$$

(c): \ Use that $\norm{u\psi-\psi}=sup_{v \in B_V} |\langle v,
u\psi \rangle -\langle v, \psi\rangle|$ and $B_V$ is
$G$-invariant.

(d): \ Take into account that $\g^*: E^* \to V^*$ is also a
$G$-map \wrt the adjoint actions of $G$ on $E^*$ and $V^*$. By the
definition of the adjoint map $\g^*$ for every $(v,u) \in V \times
G$ we have
$$<v,u\gamma^*(\psi)>-<v,\g^*(\psi)>=<v,\gamma^*(u \psi)>-<v,\g^*(\psi)>=
<\g(v), u\psi>-<\g(v),\psi>.$$ This equality implies that
$\gamma^{-1}(M) \subset V$ is SUC-small at $\gamma^*(\psi) \in
V^*$ (using the assumption that $M \subset E$ is a SUC-small set
at $\psi \in E^*$).
\end{proof}

%
%
%
%

\sk

For every $f \in RUC(G)$ denote by $\Acal_f$ the smallest closed
\emph{unital} (that is, containing the constants) $G$-invariant
subalgebra of $RUC(G)$ which contains $f$. Denote by $X_f$ the
Gelfand space of the algebra ${\mathcal A_f}$. We call ${\mathcal
A_f}$ and $X_f$ the \emph{cyclic $G$-algebra} and \emph{cyclic
$G$-system} of $f$, respectively (see for example, \cite[Ch. IV,
Section 5]{Vr} or \cite[Section 2]{GM1}). The corresponding
compactification $\a_f: G \to X_f$ is a $G$-compactification. That
is, the compact space $X_f$ is a left $G$-space such that $\a_f$
is a $G$-map
and the $G$-orbit of the point $\a_f(e)$ (where $e$ is the
identity of $G$) is dense in $X_f$. Since $f \in \Acal_f$ there
exists a continuous function $F: X_f \to \R$ such that the
following diagram commutes:
\begin{equation*}
\xymatrix { G \ar[dr]_{f} \ar[r]^{\a_f} & X_f
\ar[d]^{F} \\
  & \R }
\end{equation*}

The following theorem is one of the main results of the present
paper having in our opinion its own interest.

\begin{thm} \label{t:matUC} For a topological group $G$ and a function
$f:G\to \R$ the following conditions are equivalent:
 \ben
\item
$f \in UC(G)$.
\item The function $f: G \to \R$ is a (bicontinuous) matrix coefficient
for some continuous co-birepresentation $h_1: G \to Iso(V), \ \
h_2: G \to Iso(W)$ in a regular bilinear mapping $\om: V \times W
\to \R$ with Banach spaces $V$ and $W$.
 \een
Moreover in the second claim we can always assume without
restriction of generality that $d(V) \leq d(G)$ and $d(W) \leq
d(G)$.
\end{thm}
\begin{proof} (2) $\Longrightarrow$ (1): \ Apply Lemma
\ref{introversion}.

(1) $\Longrightarrow$ (2): \ Let $f \in UC(G)$. In particular, $f
\in RUC(G)$. Consider the corresponding cyclic algebra $\Acal_f$
and the $G$-compactification $\a_f: G \to X_f$. As we already
mentioned, $F \circ \a_f =f$ for some continuous function $F: X_f
\to \R$. Denote by $z$ the image of the identity $e \in G$ in
$X_f$ under the map $\a_f$; that is, $z:=\a_f(e)$. The compact
space $X_f$ is naturally embedded into $C(X_f)^*$ by assigning to
$y \in X_f$ the corresponding point measure $\delta_y \in
C(X_f)^*$, where $\delta_y(\varphi):=\varphi(y)$ for every
$\varphi \in C(X_f)$. In the sequel we will identify $X_f$ with
its natural image in $C(X_f)^*$.

Now we show that the $G$-orbit $FG$ of the vector $F$ in $C(X_f)$,
as a family of functions, is SUC-small (see Definition
\ref{d:small}) at $z \in X_f \subset C(X_f)^*$. Indeed, let $\eps
>0$. By our assumption, $f \in UC(G)$. In particular, $f \in
LUC(G)$. Therefore there exists a \nbd $U$ of the identity $e$ in
$G$ such that
$$|f(gu) -f(g)| < \eps
 \ \ \ \forall \ (u,g) \in U \times G.
 $$
On the other hand, since $\a_f$ is a $G$-map the equality $F \circ
\a_f =f$ implies that $$|F(guz) -F(gz)|=|F(gu
\a_f(e))-F(g\a_f(e))|=|F(\a_f(gu))-F(\a_f(g))|=|f(gu) -f(g)|.$$
Therefore we get
$$|F(guz) -F(gz)|=|<Fg,uz>-<Fg,z>| < \eps
 \ \ \ \forall \ (u,g) \in U \times G.
 $$
This means that $FG$ is SUC-small at $z$.

Let $Y:=co(-FG \cup FG)$ be the convex hull of the symmetric set \
$-FG \cup FG$. Then $Y$ is a convex symmetric $G$-invariant subset
in $C(X_f)$. Denote by $E$ the Banach subspace of $C(X_f)$
generated by $Y$, that is $E$ is the norm closure of the linear
span $sp(Y)$ of $Y$ in $C(X_f)$. Since $X_f$ is a compact
$G$-space the natural right action of $G$ on $C(X_f)$ (by linear
isometries) is continuous. By our construction $E$ is a
$G$-invariant subspace. Hence the restricted action of $G$ on $E$
is well defined and also continuous.

Since $Y$ is convex and symmetric, we can apply the construction
of Davis, Figiel, Johnson and Pelczy\'nski \cite{DFJP}. We mostly
use the presentation and the development given by Fabian in the
book \cite{Fa}. Consider the sequence $K_n:=2^n Y + 2^{-n} B_E, \
\ n \in \N$ of subsets in $E$. Let $\norm{\cdot}_n$ be the
Minkowski's functional of the set $K_n$. That is,
$$\| v\|_n = \inf\ \{\lambda
> 0 \bigm| v\in \lambda K_n\}$$
 Then $\norm{\cdot}_n$ is a norm on $E$ equivalent to the given norm
of $E$ for every $n \in \N$. For $v\in E,$ let
$$N(v):=\left(\sum^\infty_{n=1} \| v \|^2_n\right)^{1/2} \hskip 0.1cm \text{and}
\hskip 0.1cm \hskip 0.1cm V: = \{ v \in E \bigm| N(v) < \infty
\}.$$ Denote by $j: V \hookrightarrow E$ the inclusion map. Then
$(V, N)$ is a Banach space and $j: V \to E$ is a continuous linear
injection. Furthermore we have
$$F \in Y \subset j(B_V) $$
Indeed, if $y \in Y$ then $2^ny \in K_n$. So, $\| y\|_n \leq
2^{-n}$. Thus, $N(y)^2 \leq \sum_{n \in \N} 2^{-2n} <1$.

\vskip 0.3cm

By our construction $Y$ and $B_E$ are $G$-invariant. This implies
that the natural right action $V \times G \to V, \ \ (v,g) \mapsto
vg$ is isometric, that is $N(vg)=N(v)$. Moreover, by the
definition of the norm $N$ on $V$ (use the fact that the norm
$\norm{\cdot}_n$ on $E$ is equivalent to the given norm of $E$ for
each $n \in \N$) we can show that this action is norm continuous.
Therefore, the continuous co-representation $h_1: G \to Iso(V), \
h_1(g)(v):=vg$ on the Banach space $(V,N)$ is well defined.


Let $j^*: E^* \to V^*$ be the adjoint map of $j: V \to E$ and
$i^*: C(X_f)^* \to E^*$ be the adjoint of the inclusion $i: E
\hookrightarrow C(X_f)$. The composition $\gamma^*=j^* \circ i^*:
C(X_f)^* \to V^*$ is the adjoint of $\gamma:=i \circ j: V \to
C(X_f)$. Denote by $\psi$ the vector $\gamma^*(z) \in
\gamma^*(X_f) \subset V^*$. Now our aim is to show the
$G$-continuity of the vector $\psi \in V^*$, that is the
continuity of the orbit map $\widetilde{\psi}: G \to V^*$.

\sk

\nt \textbf{Claim} : $j(B_V) \subset \bigcap_{n \in \N} K_n =
\bigcap_{n \in \N} (2^n Y + 2^{-n}B_E)$.
\begin{proof}
The norms $\norm{\cdot}_n$ on $E$ are equivalent to each other. It
follows that if $v \in B_V$ then $\| v\|_n < 1$ for all $n \in
\N$. That is, $v \in \lambda_n K_n$ for some $0 < \lambda_n <1$
and $n \in \N$. By the construction $K_n$ is a convex subset
containing the origin. This implies that $\lambda_n K_n \subset
K_n$. Hence $v \in K_n$ for every $n \in\N$.
\end{proof}

Recall now that $FG$ is SUC-small at $z \in C(X_f)^*$. By Lemma
\ref{l:small}(a) we know that then $Y:=co(-FG \cup FG)$ is also
SUC-small at $z \in C(X_f)^*$. Moreover by Lemma \ref{l:small}(b)
we obtain that $M:=\bigcap_{n \in \N} (2^n Y + 2^{-n}B_E) \subset
C(X_f)$ is SUC-small at $z \in C(X_f)^*$. The linear continuous
operator $\gamma: V \to C(X_f)$ is a $G$-map. Then by Lemma
\ref{l:small}(d) it follows that $\gamma^{-1}(M) \subset V$ is
SUC-small at $\psi:=\gamma^*(z) \in V^*$. The same is true for
$B_V$ because by the claim we have $\gamma(B_V)=j(B_V) \subset M$
(and hence, $B_V \subset \g^{-1}(M)$). Now Lemma \ref{l:small}(c)
means that the orbit map $\widetilde{\psi}: G \to V^*$ is
$G$-continuous.

Define $W$ as the Banach subspace of $V^*$ generated by the orbit
$G\psi$ in $V^*$. More precisely, $W$ is the norm closure
$cl(sp(G\psi))$ of the linear span $sp(G\psi)$ of $G\psi$ in
$V^*$. Clearly, $W$ is a $G$-invariant subset of $V^*$ under the
adjoint action of $G$ on $V^*$. The left action of $G$ on $W$ by
linear isometries defines the representation $h_2: G \to Iso(W)$.
Moreover, since $\psi$ is $G$-continuous it is easy to see that in
fact every vector $w \in W$ is $G$-continuous. This means that
$h_2$ is continuous. Define the bilinear mapping $\om: V \times W
\to \R$ as a restriction of the canonical form $V \times V^* \to
\R$. Clearly, $\om$ is regular (hence, continuous) and the pair
$(h_1,h_2)$ is a continuous co-birepresentation of $G$ in $\om$.


By our construction $F \in j(V)$ (because $F \in Y \subset j(B_V)
$). Since $j$ is injective the element $v:=j^{-1}(F)$ is uniquely
determined in $V$. We already proved that $\psi=\gamma^*(z) \in
V^*$ is a $G$-continuous vector. In order to complete the proof it
suffices to show that $f=m_{v,\psi}$. Using the equality $F \circ
\a_f =f$ and the fact that $\a_f$ is a $G$-map we get $$<Fg,
z>=F(g\a_f(e))=(F \circ \a_f)(g)=f(g).$$ On the other hand,
$$
m_{v,\psi}(g)=<vg,\psi>=<j^{-1}(F)g,\gamma^*(z)>=<\gamma(j^{-1}(F)g),
z>=<Fg, z>.
$$
Hence, $f=m_{v,\psi}$, as required.

This proves the equivalence (1) $\Longleftrightarrow (2)$.
\ssk

By the $G$-continuity of $\psi$ in $W=cl(sp (G\psi))$ we get that
$d(W)=d(G\psi) \leq d(G)$. Now we check that $d(V) \leq d(G)$.
First of all $E$ is a subspace of $C(X_f)$ generated by $Y:=co(-FG
\cup FG)$. Since $F$ is a $G$-continuous vector in $E$ we have
$d(FG) \leq d(G)$. Therefore we get that
$$d(E)=d(Y)=d(FG) \leq d(G).$$
 Now it suffices to show that $d(V)=d(E)$. That is we have to show
 that the canonical construction of \cite{DFJP} in
fact always preserves the density.
Indeed, by the construction $V$ is a (diagonal) subspace into the
$l_2$-sum $Z:=\sum_{n=1}^{\infty} ((E,\norm{\cdot}_n)_{l_2}$. So,
$d(V) \leq d(Z)$. On the other hand we know that every norm
$\norm{\cdot}_n$ is equivalent to the original norm on $V$. Hence,
$d(E,\norm{\cdot}_n) \leq d(G)$. Therefore, $Z$ is an $l_2$-sum of
countably many Banach spaces each of them having the density
$d(E)$. It follows that $$d(V) \leq d(Z) \leq d(E).$$ So we obtain
$d(V) \leq d(G)$, as required.
\end{proof}


\sk
\section{Additional properties of strong dualities}
\sk

In this section we give some additional auxiliary technical
results about strong dualities (they did not appear in
\cite{me-min}). The meaning of Theorem \ref{t:crucial} is that
every continuous co-birepresentation in a bilinear mapping
naturally leads to a continuous co-birepresentation into strong
duality preserving the matrix coefficients.

\sk

\begin{defn} \label{d:major}
Let $\omega: V \times W \to \R, \ \text{and} \ \omega': V' \times
W' \to \R $ be two continuous bilinear mappings defined for
(semi)normed spaces. \ben
\item
 We say that $\om'$ \emph{refines} $\om$
(notation: $\om \succeq \om'$) if there exist continuous linear
operators of normed spaces $p: V \to V', \ q: W \to W'$ such that
$\om'(p(v),q(f))=\om(v,f)$ for every $(v,f) \in V \times W$. Hence
the following diagram commutes
\begin{equation*}
\xymatrix{ V \ar@<-2ex>[d]_{p} \times W
\ar@<2ex>[d]^{q} \ar[r]^{\om}  & \R \ar[d]^{1_{\R}} \\
V' \times W' \ar[r]^{\om'}  &  \R }
\end{equation*}
\item
Let $\Psi=(h_1,h_2)$ and $\Psi'=(h_1',h_2')$ be two
co-birepresentations of $G$ into $\om$ and $\om'$ respectively. We
say that $\Psi'$ \emph{refines} $\Psi$ (notation: $\Psi \succeq
\Psi'$) if $\om \succeq \om'$ and one can find $p$ and $q$
satisfying the assumption of the first definition such that $p$
and $q$ are $G$-maps.
\item
If $p$ and $q$ are onto then we say that $\om'$ is an \emph{onto
refinement} of $\om$. \emph{Dense refinement} will mean that
$p(V)$ and $q(W)$ are dense in $V'$ and $W'$ respectively.
Similarly, we define \emph{onto refinement} and \emph{dense
refinement} of co-birepresentations. \een
\end{defn}

\begin{lem} \label{l:sameMAT}
\ben
\item
 If $\Psi \succeq \Psi'$ then ${\mathcal
M}_{\Psi} \subset {\mathcal M}_{\Psi'}$.
\item
If $\om \succeq \om'$ is a dense refinement and $\om$ is left
(right) strong duality then $\om'$ is also left (resp., right)
strong duality.
 \een
\end{lem}
\begin{proof} (1): \
For every pair $(v, \psi) \in V \times W$ and $g \in G$ we have
$$m_{v,\psi}(g)=\om(vg,\psi)=\om'(p(vg), q(\psi))=\om'(p(v)g, q(\psi))=m_{p(v),q(\psi)}(g).$$
Thus, $m_{v,\psi}=m_{p(v),q(\psi)}$. This proves the inclusion
${\mathcal M}_{\Psi} \subset {\mathcal M}_{\Psi'}$.

(2): \ Assume that $\om$ is left strong. We show that then $\om'$
is also left strong (we omit the similar details for the "right
strong case"). Let $v_n$ be an unbounded sequence in $V'$. Since
$p(V)$ is dense in $V'$ there exists a sequence $x_n$ in $V$ such
that $\norm{p(x_n)-v_n} \leq 1$. Clearly $x_n$ is also unbounded
(otherwise $v_n$ is bounded) because $p$ is a bounded operator. By
the continuity of $q: W \to W'$ we can choose $\eps >0$ such that
$\norm{q(f)} \leq 1$ whenever $\norm{f} \leq \eps$. Since $\om$ is
left strong then $\{f(x_n)=\om(x_n,f): \ n \in \N, \ \norm{f} \leq
\eps \}$ is unbounded in $\R$. By the inclusion
$$\{f(x_n): \ n \in \N, \ \norm{f} \leq \eps \} =
\{\om'(p(x_n),q(f)): \ n \in \N,
 \ \norm{f} \leq \eps \} \subset$$
$$ \subset \{\om'(p(x_n), \phi) : \ n \in \N, \ \phi \in B_{W'}\}$$
the set $\{\phi(p(x_n)) : \ n \in \N, \ \phi \in B_{W'}\}$ is
unbounded, too. By the continuity of $\om'$ and Lemma \ref{l:cont}
there exists a constant $c >0$ such that
$$
\norm{\phi(v_n)-\phi(p(x_n)} \leq c \cdot \norm{\phi} \cdot
\norm{v_n-p(x_n)} \leq c
$$
holds for every $\phi \in B_{W'}$ and $n \in \N$. It follows that
the set $\{\phi(v_n) : \ n \in \N, \ \phi \in B_{W'}\}$ is also
unbounded in $\R$. This means that $\om'$ is also left strong.
\end{proof}

\sk

\begin{thm} \label{t:crucial} For every continuous co-birepresentation $\Psi$
of $G$ into a continuous bilinear mapping  $$\omega: V \times W
\to \R, \quad (v,f) \mapsto \om(v,f)=<v,f>=f(v)$$ such that $V$
and $W$ are normed spaces there exists a regular strong duality
$\om_0: V_0 \times W_0 \to \R$ with normed spaces $V_0$ and $W_0$
and a continuous co-birepresentation $\Psi_0$ of $G$ into $\om_0$
such that $\Psi \succeq \Psi_0$ is an onto refinement.
\end{thm}
\begin{proof}
Define a seminorm $\norm{\cdot}_*$ on $V$ by
$$
\norm{v}_*:=sup\{<v,f> : \ \ f \in B_W\}.
$$
Note that the seminorm $\norm{\cdot}_*$ on $V$ in fact is the
\emph{strong polar topology} $\beta(V,W)$ (see for
example, \cite[Section 9.4]{NB}) induced on $V$ by the form $\om:
V \times W \to \R$.

\ssk

\nt \textbf{Assertion 1}: \ \emph{There exists a constant $c >0$
such that $\norm{v}_* \leq c \cdot \norm{v}$ for every $v \in V$.}
\begin{proof}
\ By Lemma \ref{l:cont} we have the inequality $|<v,f>| \leq c
\cdot \norm{v} \cdot \norm{f}$ for some constant $c >0$. Then $<v,
\frac{f}{\norm{f}}> \leq c \cdot \norm{v}$ for every $f \in W$.
Hence, $$\norm{v}_*:=sup\{<v,\varphi> : \ \ \norm{\varphi}=1\} \
\leq c \cdot \norm{v} \  \ \quad \forall \ v \in V.$$
\end{proof}

 \ssk

\nt \textbf{Assertion 2}: \ $\om_*: (V,\norm{\cdot}_*) \times W
\to \R, \ (v,f) \mapsto \om(v,f)$ is a regular (continuous)
bilinear map.
\begin{proof}
Clearly, $<v, \frac{f}{\norm{f}}> \leq \norm{v}_*$ for every $f
\in W$ and $v \in V$. So, $<v, f> \leq \norm{v}_* \cdot \norm{f}$.
Thus, $\om_*: (V,\norm{\cdot}_*) \times W \to \R$ is regular
(hence, continuous, by Lemma \ref{l:cont}).
\end{proof}

\ssk

\nt \textbf{Assertion 3}: \ $\om_*: (V,\norm{\cdot}_*) \times W
\to \R$ is a left strong duality and the pair of natural identity
maps $(1_V)_*: V \to (V,\norm{\cdot}_*)$, $1_W: W \to W$ defines
the natural onto refinement $\om \succeq \om_*$.
\begin{proof}
Let $v_n$ be a norm unbounded sequence in $(V, \norm{\cdot}_*)$.
Then by the definition of the seminorm $\norm{\cdot}_*$ for every
$n \in \N$ there exists $f_n$ in the unit ball $B_W$ such that the
sequence $f_n(v_n)$ is unbounded in $\R$. This proves that $\om_*$
is a left strong duality.

By Assertion 1 the new $\norm{\cdot}_*$-seminorm
 topology on $V$ is coarser than the original norm
topology. It follows that the pair $(1_V)_*: V \to
(V,\norm{\cdot}_*)$, $1_W: W \to W$ defines the natural onto
refinement $\om \succeq \om_*$.
\end{proof}

\ssk

For the seminormed space $(V,\norm{\cdot}_*)$ denote by $(V_0,
\norm{\cdot}_0)$ the corresponding universal normed space. The
elements of $V_0$ can be treated as the subsets $[v]:=\{x \in V: \
\norm{x-v}_*=0 \}$ of $V$, where $v \in V$. The canonical norm on
$V_0$ is defined by $\norm{[v]}_0:=\norm{v}_*$. Denote by
$\lambda_*: (V, \norm{\cdot}_*) \to V_0, \ v \mapsto [v]$ and
$\lambda: V \to V_0, \ v \mapsto [v]$ the corresponding natural
linear continuous onto operators.

\ssk

\nt \textbf{Assertion 4}: \
The bilinear mapping
$$ \omega_L: (V_0, \norm{\cdot}_0) \times W \to \R, \ \ ([v],f)
\mapsto \om(v,f)$$ is a well defined left strong regular duality
and the pair $\lambda: V \to (V_0,\norm{\cdot})_0$, $1_W: W \to W$
defines the natural onto refinement $\om \succeq \om_L$.
Furthermore, if $\om$ is a right strong duality then $\om_L$ is a
strong duality.
\begin{proof}
Since $\R$ is Hausdorff the continuity of $\om_*$ implies that if
$\norm{v}_*=0$ then $f(v)=0$ for every $f \in W$. Hence,
$f(v_1)=f(v_2)$ for every $v_1, v_2 \in [v]$ and $f \in W$.
 This proves that
$\omega_L$ is well defined. Since
$\om_L([v],f)=\om(v,f)=\om_*(v,f)$, we easily get by Assertion 2
that $\omega_L$ is regular. Moreover the pair $(\lambda_*,1_W)$
defines the (onto, of course) refinement $\om_* \succeq \om_L$.
The latter fact implies that $\om_L$ also is left strong by Lemma
\ref{l:sameMAT}.2 and Assertion 3. Since $\om \succeq \om_*$,
$\om_* \succeq \om_L$ and $\lambda=\lambda_* \circ (1_V)_*$ we get
the natural onto refinement $\om \succeq \om_L$ with respect to
the pair $(\lambda, 1_W)$.

Now if $\om$ is a right strict duality then $\om_L$ remains right
strong duality by Lemma \ref{l:sameMAT}.2.
\end{proof}

Now assume that the pair $h_1: G \to Iso(V, \norm \cdot)$, \ $h_2:
G \to Iso(W)$ is a given continuous co-birepresentation $\Psi$ of
$G$ in $\om$. Since $\norm{v}_* \leq c \cdot \norm{v}$ we get that
every vector $v \in V$ is $G$-continuous \wrt $\norm{\cdot}_*$. On
the other hand, since $B_W$ is $G$-invariant we get that
$\norm{vg}_*=\norm{v}_*$ for every $g \in G$. Therefore, $(V,
\norm{\cdot}_*) \times G \to (V, \norm{\cdot}_*), \
(v,g):=vg=h_1(g)(v)$ is a well defined continuous right action.
%

Define the co-representation $(h_1)_0: G \to Iso(V_0,
\norm{\cdot}_0)$ by the natural right action $([v],g) \mapsto
[v]g=[vg]$. Then it is a well defined continuous
co-representation. The proof is straightforward taking into
account the trivial equalities $\norm{[v]}_0=\norm{v}_*$, $[v]g
=[vg]$ and $\norm{vg}_*=\norm{v}_*$ for every $(v,g) \in V \times
G$. It is also easy to see that the equality
$\om_L([v]g,\psi)=\om_L([v],g\psi)$ holds for every $([v],g,\psi)
\in V_0 \times G \times W$. Hence, the pair $((h_1)_0,h_2)$
defines a continuous co-birepresentation (denote it by $\Psi_L$)
of $G$ in $\om_L: (V_0, \norm{\cdot}_0) \times W \to \R$.
Furthermore, each of the maps $\lambda$ and $1_W$ from Assertion 4
are $G$-maps. So in fact we found a co-birepresentation $\Psi_L$
into the left strong regular duality $\om_L$ such that $\Psi
\succeq \Psi_L$.

\ssk

Similarly, starting now from $\om_L$ and switching left and right,
we can construct: a seminormed space $(W, \norm{\cdot}_*)$, its
universal normed space $W_0$, a regular \emph{right strong}
duality $(\om_L)_R: V_0 \times W_0 \to \R$ which is an onto
refinement of $\om_L$ and a continuous co-birepresentation
$(\Psi_L)_R$ of $G$ into $(\om_L)_R$ such that $\Psi_L \succeq
(\Psi_L)_R$. Denote by $\om_0$ the duality $(\om_L)_R$ and by
$\Psi_0$ the co-birepresentation $(\Psi_L)_R$. Then $\Psi_0$ is
the desired co-birepresentation because $\Psi \succeq \Psi_L
\succeq \Psi_0$ and $\om_0$ in fact is left and right strong (take
into account the analogue of Assertion 4).
\end{proof}

\begin{lem} \label{l:completion}
\ben
\item
 Let \ $\omega: V \times W \to \R$ be a strong duality such
that $V$ and $W$ are normed spaces. Denote by $\varrho_V: V \to
\widehat{V}$ and $\varrho_W: W \to \widehat{W}$ the corresponding
completions. Then the uniquely defined continuous extension
$$ \widehat{\omega}: \widehat{V} \times \widehat{W} \to \R
$$
is a strong duality and the pair $(\varrho_V, \varrho_W)$ defines
the canonical refinement $\om \succeq \widehat{\omega}$. If $\om$
is regular then the same is true for $\widehat{\omega}$.
\item
Assume that the pair $h_1: G \to Iso(V)$, \ $h_2: G \to Iso(W)$
defines a continuous co-birepresentation $\Phi$. Then there exists
a uniquely defined extension to a continuous co-birepresentation
$\widehat{\Phi}$ of $G$ into the form $\widehat{\omega}$ defined
by the pair $\widehat{h_1}: G \to Iso(\widehat{V})$, \
$\widehat{h_2}: G \to Iso(\widehat{W})$. In fact, $\widehat{\Phi}$
is a dense refinement of $\Phi$.
 \een
\end{lem}
\begin{proof} (1) It can be derived by Lemma \ref{l:sameMAT}.2
because $\widehat{\omega}$ is a dense refinement of $\om$ under
the natural dense inclusions $\varrho_V: V \to \widehat{V}$ and
$\varrho_W: W \to \widehat{W}$. The regularity of
$\widehat{\omega}$ for regular $\om$ is clear.


(2) \ Note that for every continuous (not necessarily isometric)
linear (left or right) action of a topological group $G$ on a
normed space $V$ there exists a uniquely defined canonical linear
extension on the Banach space $\widehat{V}$ which is also
continuous. This is easy to verify directly or it can be derived
also by \cite[Proposition 2.6.4]{me-fr}. Straightforward arguments
show also that: $\varrho_V$ and $\varrho_W$ are continuous
$G$-maps, $\widehat{\omega}(xg,y)=\widehat{\omega}(x,gy)$ for
every $(x,g,y) \in \widehat{V} \times G \times \widehat{W}$ and
the corresponding $g$-translations $\widehat{V} \to \widehat{V}$
and $\widehat{W} \to \widehat{W}$ are linear isometries.
\end{proof}

\begin{prop} \label{p:strongization}
Let $G$ be a Hausdorff topological group and let
$$\Phi:=\{\Phi_i\}_{i \in I}=\{\om_i: E_i \times F_i \to \R\, \quad
\a_{1i}: E_i \times G \to E_i, \ \a_{2i}: G \times F_i \to
F_i\}_{i \in I}$$ be a system of continuous co-birepresentations
of $G$ into regular bilinear mappings $\om_i$ with Banach spaces
$E_i, \ F_i$. Let
$$
{\mathcal M}_{\Phi}:=\{m_{v,\psi}: G \to \R:  \ \ \ \ (v,\psi) \in
E_i \times F_i, \ i \in I \}
$$
be the family of all corresponding matrix coefficients. Suppose
that ${\mathcal M}_{\Phi}$ is a local separating family of
functions on $G$. Then there exists a continuous t-exact
co-birepresentation $\Psi$ of $G$ into a regular strong duality $
\omega: E \times F \to \R$. Furthermore we can assume that $E$ and
$F$ are Banach spaces and their densities are not greater than
$sup\{d(E_i) \cdot d(F_i) \cdot |I|: \ \ i \in I \}$.
\end{prop}
\begin{proof}
Consider the $l_2$-sum of the given system $\Phi$ of
co-birepresentations. That is, define naturally the Banach spaces
$V:=(\sum_{i \in I} E_i)_{l_2}$ and $W:=\sum_{i \in I}
(F_i)_{l_2}$, the continuous co-representation $h_1: G \to Iso(V)$
and the continuous representation $h_2: G \to Iso(W)$. Clearly,
$<vg,f>=<v,gf>$ for the natural bilinear mapping
$$\om_{l_2}: V \times W \to \R, \ \ <\sum_i v_i, \sum_i
f_i>:=\sum_if_i(v_i).$$

Since $|\om_i(v,f)| \leq \norm{v} \cdot \norm{f}$ for every $i \in
I$ it follows by the Schwartz inequality that the form $\om_{l_2}$
is well defined and continuous. Moreover the co-birepresentation
$\Psi_{l_2}:=(h_1,h_2)$ of $G$ in $\om_{l_2}$ is also well defined
and continuous.

Now in order to get a regular strong duality we apply Theorem
\ref{t:crucial} to $\Psi_{l_2}$. Then we obtain the regular strong
duality $(\om_{l_2})_0: V_0 \times W_0 \to \R$ and a continuous
co-birepresentation $(\Psi_{l_2})_0$
of $G$ in $(\om_{l_2})_0$ such that $\om_{l_2} \succeq
(\om_{l_2})_0$ and $\Psi_{l_2} \succeq (\Psi_{l_2})_0$. Applying
Lemma \ref{l:completion} we get the continuous co-birepresentation
$\widehat{(\Psi_{l_2})_0}$ of $G$ into a regular strong duality
$\widehat{(\om_{l_2})_0}: \widehat{V_0} \times \widehat{W_0} \to
\R$ such that $(\om_{l_2})_0 \succeq \widehat{(\om_{l_2})_0}$ and
$(\Psi_{l_2})_0 \succeq \widehat{(\Psi_{l_2})_0}$. We claim that
$\Psi:= \widehat{(\Psi_{l_2})_0}$ is the desired
co-birepresentation into $\om:=\widehat{(\om_{l_2})_0}:
\widehat{V_0} \times \widehat{W_0} \to \R$. Indeed, first of all
observe that for every $i \in I$ the co-representation $\Phi_i$ of
$G$ in $\om_i$ can be treated as "a part" of the global
co-birepresentation $\Psi_{l_2}$. Therefore the set ${\mathcal
M}_{\Psi_{l_2}}$ of all matrix coefficients defined by
$\Psi_{l_2}$ contains the set ${\mathcal M}_{\Phi}$ which is a
local separating family on $G$. Hence ${\mathcal M}_{\Psi_{l_2}}$
is also local separating. Now by Lemma \ref{l:sameMAT}.1 the same
is true for the families ${\mathcal M}_{(\Psi_{l_2})_0}$ and
${\mathcal M}_{\Psi} = {\mathcal M}_{\widehat{(\Psi_{l_2})_0}}$
because $\Psi_{l_2} \succeq (\Psi_{l_2})_0 \succeq
\widehat{(\Psi_{l_2})_0}=\Psi$. It follows by Lemma \ref{l:Mat}
that the co-birepresentation $\Psi$ of $G$ in
$\om=\widehat{(\om_{l_2})_0}: E \times F \to \R$ is t-exact, where
$E:=\widehat{V_0}$ and $F:=\widehat{W_0}$ are certainly Banach
spaces. The completion of normed spaces does not increase the
density. So by our construction (using some obvious properties of
$l_2$-sums) one can assume in addition that the densities of $E$
and $F$ are not greater than $sup\{d(E_i) \cdot d(F_i) \cdot |I|:
\ \ i \in I \}$.
\end{proof}

 \sk
\section{Proof of the main theorem and some consequences}
\label{s:mainthm}
 \sk

First we prove the following crucial result.

\begin{thm} \label{t:allBR}
\ben \item
 Every Hausdorff topological group $G$ is a BR-group.
\item
 Moreover,
there exists a t-exact birepresentation
$$\Psi:=\{\om: E \times F \to \R, \ h_1: G \to
Iso(E), \ h_2: G \to Iso(F)\}
$$
of $G$ such that: $\om$ is a regular strong duality; $E$ and $F$
are Banach spaces with the density not greater than $w(G)$. \een
\end{thm}
\begin{proof} \ (1) directly follows from (2). Hence it suffices to show (2).
By Lemma \ref{l:Roelcke} the algebra $UC(G)$ separates points and
closed subsets of $G$. Choose a subfamily $S \subset UC(G)$ with
cardinality $|S| \leq \chi(G)$ such that $S$ is a local separating
family (see Definition \ref{d:local-separ}) for $G$. By Theorem
\ref{t:matUC} every $f \in S \subset UC(G)$ can be represented as
a bicontinuous matrix coefficient by some \emph{regular} bilinear
mapping. More precisely, there exists a continuous
co-birepresentation $\Phi_f$ defined by the pair $h_f: G \to
Iso(V_f), \ \ h_f': G \to Iso(W_f)$ in a regular bilinear mapping
$\om_f: V_f \times W_f \to \R$ with Banach spaces $V_f$ and $W_f$
such that $f=m_{v,\psi}$ for some pair $(v,\psi) \in V_f \times
W_f$. Moreover we can assume that $d(V_f) \leq d(G)$ and $d(W_f)
\leq d(G)$.

We get a system
$$\Phi:=\{\Phi_f \}_{f \in S}=\{\om_f: V_f
\times W_f \to \R,  \quad h_f: G \to Iso(V_f), \ \ h'_f: G \to
Iso(W_f) \}_{f \in S}
$$
of continuous co-birepresentations of $G$.
 By our construction the
corresponding set of all matrix coefficients ${\mathcal M}_{\Phi}$
contains the local separating family $S$ of functions on $G$. We
can apply Proposition \ref{p:strongization}. Then there exists a
continuous t-exact co-birepresentation $\Psi'$ of $G$ into a
regular \emph{strong} duality $ \omega: E \times F \to \R$ with
Banach spaces $E$ and $F$ the densities of them are not greater
than $sup\{d(V_f) \cdot d(W_f) \cdot |S|: \ \ f \in S \}$. Since
$|S| \leq \chi(G)$, $w(G)=d(G) \cdot \chi(G)$, and $d(V_f) \leq
d(G)$, $d(W_f) \leq d(G)$, we obtain that $max\{d(E),d(F)\} \leq
w(G)$. Finally, in order to get the desired
\emph{birepresentation} $\Psi$ from our co-birepresentation
$\Psi'$, just define it according to Remark \ref{r:opp} as
$\Psi:=(\Psi')^{op}$.
\end{proof}

Now we obtain our main result:

\begin{thm} \label{t:main}
 For every Hausdorff topological group $G$ there exists a minimal group
$M$ and a topological group retraction $p: M \to G$.

Furthermore we can assume that: \bit
\item [(a)]
There exists a regular strong duality $\om: E \times F \to \R$
with Banach spaces $E, \ F$ and a t-exact birepresentation $\Psi$
of $G$ into $\om$ such that: $max\{d(E),d(F)\} \leq w(G)$ and $M$
can be constructed as the induced group $M_+(\Psi)$ of $\Psi$.
That is
$$M:=M_+(\Psi)=H_+(\om) \leftthreetimes_{\pi} G=((\R \times E \leftthreetimes_{\om^{\nabla}} F)
\leftthreetimes_{\a} \R_+) \leftthreetimes_{\pi} G.$$
\item [(b)]
$p: M \to G$ is the natural group retraction and the group
$H_+(\om)=ker(p)$ is minimal, Weil complete and solvable.
\item [(c)]
 $w(M)=w(G)$, $\chi(M)=\chi(G)$ and
$\psi(M)=\psi(G)$.
\item [(d)]
If $G$ is Raikov complete (Weil complete) then $M$ also has the
same property.
\item [(e)]
If $G$ is solvable then $M$ is solvable. \eit
\end{thm}
\begin{proof} \ (a): \
Use Theorems \ref{t:allBR} and \ref{t:BR} (see also Theorem
\ref{t:BRapp} below).

(b): \ By Proposition \ref{str-dual-ismin} the kernel
$ker(p)=H_+(\om)$ of the retraction $p$ is a minimal group.
Furthermore, using Lemma \ref{l:complconcrete} we can conclude
that the group $H_+(\om)$ is Weil complete. The solvability of the
group $H_+(\om)=(\R \times E \leftthreetimes_{\om^{\nabla}} F)
\leftthreetimes_{\a} \R_+$ is trivial.

(c): \ Take into account that $M$, as a topological space, is
homeomorphic to the product of $G$ and the metrizable space
$H_+(\om)$ with weight $\leq w(G)$.

(d): \ Use Lemmas \ref{l:completionofsemid} and
\ref{l:complconcrete}.

(e): \ Use (b).
\end{proof}

\sk

In particular, by Theorem \ref{t:main} we can conclude now that
Pestov's conjecture is true.

\sk \ssk

\begin{cor} \label{ex:arh-pr}
Every compact Hausdorff homogeneous space admits a transitive
continuous action of a minimal group.
\end{cor}
\begin{proof} \
Let $X$ be a compact Hausdorff homogeneous space. Then the group
$G:=Homeo(X)$ of all homeomorphisms of $X$ is a topological group
\wrt the usual compact open topology. The natural action $\a: G
\times X \to X$ is continuous. This action is transitive because
$X$ is homogeneous. By Theorem \ref{t:main} there exists a minimal
group $M$ and a continuous group retraction $p: M \to G$. Then the
action $M \times X \to X, \ mx:=\a(p(m),x)$ is also continuous and
transitive.
\end{proof}

\begin{remarks} \label{r:H_+}
\ben
\item
Theorem \ref{t:main} implies that there exist Raikov-complete
(Weil-complete) minimal topological groups such that $\chi(M)$ and
$\psi(M)$ are different (as far as possible for general groups).
This answers negatively Question C (see Introduction).
\item
Applying Corollary \ref{ex:arh-pr} to a not dyadic compact
homogeneous space $X$ we get an immediate negative answer to
Question D (see Introduction).
 \een
\end{remarks}


Recall that a topological group $K$ is \emph{perfectly minimal} in
the sense of Stoyanov (see for example \cite{DPS}) if the product
$K \times P$ is minimal for every minimal group $P$. By the test
of perfect minimality \cite[Theorem 1.14]{me-min} a minimal group
$K$ is perfectly minimal iff its center is perfectly minimal. It
is easy to see that the center of the group $M=H_+(\om)
\leftthreetimes_{\pi} G$ is trivial. Indeed, the center of its
subgroup $H_+(\om)$ is already trivial. Here it is important to
note that the bilinear mapping $\om$ (being a strong duality) is
separated (see the text after Definition \ref{d:str-duality}).
Therefore, it follows that in fact, $M$ in the main Theorem is
perfectly minimal.

\sk Recall that every \lca group $G$ is a group retract of a
generalized Heisenberg group $H(\Delta)=\T \times G^*
\leftthreetimes_{\Delta} G$ (see Section 2 or \cite{me-min}) which
certainly is locally compact. For nonabelian case the following
question (recorded also in \cite[Question 3.3.5]{CHR}) seems to be
open.

\begin{question} \cite[Question 2.13.1]{me-min}
Is it true that every locally compact Hausdorff group is a group
retract (quotient, or a subgroup) of a locally compact minimal
group ?
\end{question}

\ssk

Theorem \ref{t:allBR} shows that every topological group $G$
admits sufficiently many continuous representations into
continuous bilinear mappings. It turns out that in general we
cannot replace general bilinear mappings $E \times F \to \R$ by
the \emph{canonical} bilinear mappings $<,>: V \times V^* \to \R$.
More precisely, let $h: G \to Iso(V)$ be a given continuous
representation of $G$ on $V$. Then we have the adjoint
representation $h^*: G \to Iso(V^*), \ h^*(g)(\psi)=g\psi$, where
$(g\psi)(v)=\psi(g^{-1}v)$. One attractive previous idea to prove
that every group $G$ is a BR-group was trying to find sufficiently
many representations $h$ of $G$ such that $h^*$ is also
continuous. If $h^*$ is continuous then $(h,h^*)$ becomes a
continuous birepresentation in $<,>$. Every topological group can
be treated as a topological subgroup of $Iso(V)$ for a suitable
Banach space $V$ (see for example, \cite{Te,pest-wh}) the
t-exactness is clear by Remark \ref{r:definitionsand}.5. Thus, we
could derive directly that every topological group $G$ is a
BR-group. Although this result is really true (Theorem
\ref{t:allBR}) however in its proof we cannot use that direct
naive argument. The reason is that in general $h^*$ is not
continuous (see for example \cite{me-fr}).

This remark suggests the following definition.

\begin{defn} \label{d:adj}
We say that a topological group $G$ is \emph{adjoint continuously
representable} (in short: \emph{ACR-group}) if there exists a
continuous representation $h: G \to Iso(V)$ on a Banach space $V$
such that the adjoint representation $h^*: G \to Iso(V^*)$ is also
continuous and the continuous birepresentation $(h,h^*)$ of $G$ is
t-exact.
\end{defn}

It seems to be interesting to study the class of adjoint
continuous representable groups. Note that $h^*$ is continuous for
every continuous representation $h: G \to Iso(V)$ on an Asplund
(e.g., reflexive) Banach space $V$ (see \cite{me-fr}). It follows
that every Asplund representable, e.g., reflexively representable,
group $G$ is an ACR-group (where Asplund (resp., reflexively)
representability means that $G$ can be embedded into $Iso(V)$ for
some Asplund (resp., reflexive) space $V$). For instance every
locally compact Hausdorff group $G$ (being Hilbert representable)
is an ACR-group. Only recently became clear that this result
cannot cover all groups. Indeed the group $Homeo_+[0,1]$ (the
topological group of all orientation preserving homeomorphisms of
$[0,1]$) is not reflexively representable \cite{me-ref} and even
not Asplund representable \cite{GM2}. Moreover the following
stronger result \cite{GM2} is true (note that this proved also by
V. Uspenskij): if $G:=Homeo_+[0,1]$  and $h: G \to Iso(V)$ is a
continuous representation on a Banach space $V$ such that the
adjoint representation $h^*: G \to Iso(V^*)$ is also continuous
then $h$ is trivial. It follows that $Homeo_+[0,1]$ is not an
ACR-group. For more information and questions about group
representations on Banach spaces we refer to \cite{GM1, GM2,
Me-opit2}.

\sk \sk
\section{Appendix}
\label{s:appendix}
\sk

For the readers convenience we include here (sometimes simplified)
proofs of some principal results from \cite{me-min} which we need
in the present paper.

\sk

For every topological group $(P,\g)$ and its subgroup $H$ denote
by $\g/H$ the usual quotient topology on the coset space $P/H$.
More precisely, if $pr: P \to P/H$ is the canonical projection
then $\g/H:=\{OH: \ O \in \g \}=\{pr(O): \ O \in \g \}$. If $q: P
\to G$ is an onto homomorphism. Then on $G$ we can define the
quotient (group) topology which in fact is the topology $q(\g)$.

The following well known result is very useful.

\begin{lem} \label{merson} \emph{(Merson's Lemma)} Let $(G, \g)$ be a
not necessarily Hausdorff topological group and $H$ be a not
necessarily closed subgroup of $G$. Assume that $\g_1 \subset \g$
be a coarser group topology on $G$ such that $\g_1|_H=\g|_H$ and
$\g_1/H=\g /H$. Then $\g_1=\g$.
\end{lem}
\begin{proof}
See for example \cite[Lemma 7.2.3]{DPS}.
\end{proof}


\begin{defn} \label{d:2.4}
Let $q: X \to Y$ be a (not necessarily group \footnote{That is,
$q$ is not necessarily a homomorphism}) retraction of a group $X$
on a subgroup $Y$. We say that $q$ is \emph{central} if
$$
q(xyx^{-1})=y \ \ \ \forall \ (x,y) \in X \times Y.
$$
\end{defn}

\begin{lem} \label{l:2.5} Let $H(\om)=A \times E \leftthreetimes_{\om^{\nabla}} F$
be the Heisenberg group of the biadditive mapping $\om: E \times F
\to A$. Then the natural projections $q_E: H(\om) \to E, \ q_F:
H(\om) \to F$ and $q_A: H(\om) \to A$ are central.
\end{lem}
\begin{proof}
If $u:=(a,x,f) \in H(\om), \ y \in E$, $\varphi \in F$ and $a \in
A$ then $uyu^{-1}=(f(y),y,0_F)$,
$u\varphi^{-1}u^{-1}=(-\varphi(x),0_E,\varphi)$ and $uau^{-1}=a$.
\end{proof}

\begin{prop} \label{p:2.6}
Let $(M,\g)$ be a topological group such that $M$ is
\emph{algebraically} \footnote{That is the topology on $M$ is not
necessarily the product topology of $X \times G$} a semidirect
product $M = X \leftthreetimes_{\a} G$. If $q: X \to Y$ is a
continuous central retraction of the topological subgroup $X$ on a
topological $G$-subgroup $Y$ of $X$, then the action
\begin{equation} \label{f1}
\a|_{G\times Y}: (G,\g/X) \times (Y, \g|_Y) \to (Y, \g|_Y)
\end{equation}
is continuous.
\end{prop}
\begin{proof}
By our assumption $M$ \emph{algebraically} is the semidirect
product $M = X \leftthreetimes_{\a} G$. Therefore we have $M/X=\{X
\times \{g\}\}_{g \in G}$. We sometimes identify $M/X$ and $G$.
This justifies also the notation $(G,\g/X)$. Note also that
%
then the group topologies $\g/X$ and $pr(\g)$ are the same on $G$,
where $pr: M \to G=M/X, \ (x,g) \mapsto g$ denotes the canonical
projection.

Clearly, each $g$-transition $(Y, \g|_Y) \to (Y, \g|_Y)$ is
continuous. Hence it suffices to show that
the action \ref{f1} is continuous at $(e_G,y)$ for every $y \in
Y$, where $e_G$ is the neutral element of $G$. Fix an arbitrary $y
\in Y$ and a \nbd $O(y)$ of $y$ in $(Y, \g|_Y)$. Since the
retraction $q: (X,\g|_X \to (Y, \g|_Y)$ is continuous (at $y$)
there exists a \nbd $U_1$ of $y:=(y,e_G)$ in $(M,\g)$ such that
$$
q(U_1 \cap X) \subset O.
$$

The conjugation $M \times M \to M, \ (a,b) \to aba^{-1}$ is
continuous (at $(e_M,y)$). We can choose: a \nbd $U_2$ of $y$ in
$M$ and a \nbd $V$ of $e_M$ in $M$ such that
$$
vU_2v^{-1} \subset U_1  \ \ \ \ \ \ \forall \ v \in V.
$$

 Now, we claim that (for the canonical projection $pr: M \to
G=M/X$) we have
$$
\a(g,z) \in O  \ \ \ \ \ \ \forall \ z \in U_2 \cap Y, \ \ g \in
pr(V).
$$
Indeed, if $v=(x,g) \in V$ and $z \in U_2 \cap Y$ then $vzv^{-1}
\in U_1$ because $vU_2v^{-1} \subset U_1$. From the normality of
$X$ in $M$ we have $vzv^{-1} \in X$. Thus, $vzv^{-1} \in U_1 \cap
X$. Elementary computations show that
$$
vzv^{-1}=(x,g)(z,e_G)(x,g)^{-1}=(x\a(g,z)x^{-1},e_G)=x\a(g,z)x^{-1}.
$$
Using the inclusion $q(U_1 \cap X) \subset O$ we get $q(x
\a(g,z)x^{-1}) \in O$. Since $q$ is a central retraction and
$\a(g,z) \in Y$, we obtain $q(x \a(g,z)x^{-1})=\a(g,z)$.
Therefore, $\a(g,z) \in O$ for every $g \in pr(V)$ and $z \in U_2
\cap Y$. Finally observe that $pr(V)$ is a \nbd of $e_G$ in $(G,
\g/X)$ and $U_2 \cap Y$ is a \nbd of $y$ in $(Y,\g|_Y)$. This
means that the action \ref{f1} is continuous at $(e_G,y)$.
\end{proof}

\begin{prop}\label{p:2.7}
Let $M=(X \leftthreetimes_{\a} G, \g)$ be a topological semidirect
product and $\{Y_i\}_{i \in I}$ be a system of $G$-subgroups in
$X$ such that the system of actions
$$
\{\a|_{G \times Y_i}: G \times Y_i \to Y_i\}_{i \in I}\
$$
is t-exact. Suppose that for each $i \in I$ there exists a
continuous central retraction $q_i: X \to Y_i$. Then if $\g_1
\subset \g$ is a coarser group topology on $M$ such that $\g_1|_X
=\g|_X$ then $\g_1=\g$.
\end{prop}

\begin{proof}
Proposition \ref{p:2.6} shows that each action
$$
\a|_{G \times Y_i}: (G, \g_1/X) \times Y_i \to Y_i
$$
is continuous. Clearly, $\g_1/X \subset \g/X$. By Definition
\ref{d:t-exactness}.2 the group topology $\g_1/X$ coincides with
the given topology $\g/X$ of $G$. Now, Merson's Lemma \ref{merson}
implies that $\g_1=\g$.
\end{proof}

\begin{cor} \label{c:2.8}
Let $(X \leftthreetimes_{\a} G, \g)$ be a topological semidirect
product and let $\a: G \times X \to X$ be t-exact. Suppose that
$X$ is abelian and $\g_1 \subset \g$ is a coarser group topology
which agrees with $\g$ on $X$. Then $\g_1=\g$.
\end{cor}
\begin{proof}
Since $X$ is abelian, the identity mapping $X \to X$ is a central
retraction.
\end{proof}

The commutativity of $X$ is essential here as we already mentioned
in Remarks \ref{r:definitionsand}.6.

\sk

\begin{lem} \label{l:normunbounded}
Let $(E,\norm{\cdot})$ be a normed space.
Denote by $\s$ the norm topology on $E$. Suppose that $\s' \subset
\s$ is a strictly coarser, not necessarily Hausdorff, group
topology on $E$. Then every $\s'$-open nonempty subset $U$
in $E$ is norm-unbounded.
\end{lem}
\begin{proof}
Since $\s'$ is strictly coarser than the given norm topology,
there exists $\eps_0 >0$ such that every $\s'$-\nbd $O$ of $0_E$
in $E$ contains an element $x$ with $\norm{x} \geq \eps_0$. It
suffices to prove our lemma for a $\s'$-\nbd $U$ of $0_E$. Since
$\s'$ is a group topology, for each natural $n$ there exists a
$\s'$-\nbd $V_n$ such that $nV_n \subset U$. One can choose $x_n
\in V_n$ with the property $\norm{x_n} \geq \eps_0$. Then $\norm{n
\cdot x_n} \geq n \cdot \eps_0$. Since $nx_n \in U$ (and $n$ is
arbitrary), this means that $U$ is norm-unbounded.
\end{proof}

\begin{lem} \label{l:strong}
Let $\om: (E, \s) \times (F, \t) \to \R$ be a strong duality with
normed spaces $E$ and $F$. Assume that $\s' \subset \s$ and $\t'
\subset \t$ are coarser, not necessarily Hausdorff, group
topologies on $E$ and $F$ respectively such that $\om: (E, \s')
\times (F, \t') \to \R$ is continuous. Then necessarily $\s'=\s$
and $\t'=\t$.
\end{lem}
\begin{proof} We show that $\s'=\s$. We omit the similar arguments for
$\t'=\t$.

By our assumption $\om: (E, \s') \times (F, \t') \to \R$ is
continuous. Then this map remains continuous replacing $\t'$ by
the stronger topology $\t$. That is the map $\om: (E, \s') \times
(F, \t) \to \R$ is continuous, too. Assume that $\s'$ is strictly
coarser than $\s$. By Lemma \ref{l:normunbounded} every $\s'$-\nbd
of $0_E$ is norm-unbounded in $(E,\norm{\cdot})$. By the
continuity of $\om: (E, \s') \times (F, \t) \to \R$ at the point
$(0_E,0_F)$ there exist: an $\s'$-\nbd $U$ of $0_E$ and a
$\t$-\nbd $V$ of $0_F$ such that $\{f(u): \ \ u \in U, \ f \in V
\} \subset (-1,1)$. Since $U$ is norm-unbounded the set $\{f(u): \
u \in U, \ f \in V \} \subset (-1,1)$ is also unbounded in $\R$ by
Definition \ref{d:str-duality}. This contradiction completes the
proof.
\end{proof}

\begin{prop} \label{p:2.9}
Let $H(\om)=\R \times E \leftthreetimes_{\om^{\nabla}} F$ be the
Heisenberg group of the strong duality $\om: E \times F \to \R$
with normed spaces $E$ and $F$. Assume that $\g_1 \subset \g$ is a
coarser group topology on $H(\om)$ such that
$\g_1|_{\R}=\g|_{\R}$. Then $\g_1=\g$.
\end{prop}
\begin{proof}
 Denote
by $\g$ the given product topology on $H(w)$. Let $\g_1 \subset
\g$ be a coarser
group topology on $H(w)$ such that $\g_1|_{\R}=\g|_{\R}$. By
Merson's Lemma it suffices to show that $\g_1/\R=\g/\R$.

%

First we establish the continuity of the map
\begin{equation} \label{f2}
w: (E,\g_1 |_E) \times (F,\g_1 / \R \times E) \to (\R, \g_1
|_{\R})=(\R,\g|_{\R})
\end{equation}

 We proof the continuity of the map \ref{f2}
 at an arbitrary pair $(x_0, f_0) \in E \times F$. Let $O$ be a
 neighborhood of $f_0(x_0)$ in $(\R, \g_1 |_{\R})$.
 Choose a neighborhood $O'$ of $(f_0(x_0),0_E,0_F)$ in $(H(w),\g_1)$ such that
$O' \cap \R = O$.
 Consider the points
$\bar{x_0}:=(0_{\R},x_0,0_F), \bar{f_0}:=(0_{\R},0_E,f_0) \in
H(w)$. Observe that the commutator $[\bar{f_0}, \bar{x_0}]$ is
just $(f_0(x_0),0_E,0_F)$. Since $(H(w), \g_1)$ is a topological
group there exist $\g_1$-neighborhoods $U$ and $V$ of $\bar{x_0}$
and $\bar{f_0}$ respectively such that $[v,u] \in O'$ for every
pair $v \in V, u \in U$.
In particular, for every $\bar{y}:=(0_{\R},y,0_F) \in U \cap E$
and $v:=(a,x,f) \in V$ we have $[v,\bar{y}]=(f(y),0_E,0_F) \in O'
\cap \R=O$. We obtain that $f(y) \in O$ for every $f \in q_F(V)$
and $\bar{y} \in U \cap E$. This means that we have the continuity
of \ref{f2} at $(f_0,x_0)$ because $q_F(V)$ is a neighborhood of
$f_0$ in the space $(F,\g_1 / \R \times E)$ and $U \cap E$ is a
neighborhood of $x_0$ in $(E, \g_1|E)$. Since the given biadditive
mapping is a strong duality it follows by Lemma \ref{l:strong}
that the topology $\g_1 / \R \times E$ on $F$ coincides with the
given topology $\t=\g / \R \times E$.

Quite similarly one can prove that the following map is continuous
$$
w: (E,\g_1 /\R \times F) \times (F,\g_1 |_F) \to (\R, \g_1 |_{\R})
$$
Which implies that $\g_1 /{\R \times F}= \g /{\R \times F}$.

Denote by $\s$ and $\t$ the given norm topologies on $E$ and $F$
respectively.

By the equalities $\g_1 /\R \times E= \g /\R \times E= \t$ in $F$
and $\g_1 /\R \times F= \g /\R \times F=\s$ in $E$ it follows that
the maps
$$q_E: (H(w),\g_1) \to (E, \s), \hsk \ \ (a,x,f) \mapsto x$$
and
$$q_F: (H(w),\g_1) \to (F,\t), \hsk \ \ (a,x,f) \mapsto f$$
are continuous. Then we obtain
that
$$q_{E \times F}: (H(w),\g_1) \to E \times F, \hsk \ \ (a,x,f) \mapsto (x,f)$$
is also continuous, where $E \times F$ is endowed with the product
topology induced by the pair of topologies $(\s,\t)$. This
topology coincides with $\g/\R$. Then $\g_1 /\R \supset \g /\R$.
Since $\g_1 \subset \g$ we have $\g_1 /\R \subset \g /\R$. Hence
$\g_1 /\R = \g /\R$, as desired.
\end{proof}

\begin{thm} \label{t:plus} \cite[Theorem 3.9]{me-min}
For every strong duality $\om: E \times F \to \R$ with normed
spaces $E$ and $F$ the corresponding group $H_+(\om)=H(\om)
\leftthreetimes_{\a} \R_+$ is minimal.
\end{thm}
\begin{proof}
Denote by $\g$ the given topology on $H_+(\om)$ and suppose that
$\g_1 \subset \g$ is a coarser Hausdorff group topology. The group
$\R \leftthreetimes \R_+$ is minimal \cite{DS},
 (see Introduction) and it naturally is embedded in
$H_+(\om)$. Therefore, $\g_1|_{\R} = \g|_{\R}$. From Proposition
\ref{p:2.9} immediately follows that
$\g_1|_{H(\om)}=\g|_{H(\om)}$.

Now observe that by Lemma \ref{l:2.5} the natural retraction $q:
H(\om) \to \R$ is
central and the action of $\R_+$ on $\R$ is t-exact (see Remark
\ref{r:definitionsand}.2). By Proposition \ref{p:2.7} (in the
situation: $G:=\R_+, \ X:=H(\om), \ Y:=\R$) we get $\g_1=\g$.
\end{proof}

\sk

The following result is a particular case of \cite[Theorem
4.8]{me-min}. For simplicity we give the arguments only for the
system with a single birepresentation. This particular case is
enough for the main result of the present paper.

\begin{thm} \label{t:BRapp} \emph{(Compare \cite[Theorem 4.8]{me-min})}
 Let $\Phi$ be a t-exact birepresentation of a topological group $G$
 into a strong duality
 $\om: E \times F \to \R$ with normed spaces $E$ and $F$.
 \ben
\item
 Then the corresponding induced group
 $$
 M:=M_+(\Phi)=((\R \times E \leftthreetimes_{\om^{\nabla}} F)
\leftthreetimes_{\a} \R_+) \leftthreetimes_{\pi} G
$$ is minimal.
\item
The projection $p: M \to G$ is a group retraction such that $M$
and also the kernel $ker(p)$ are minimal groups.
\een
\end{thm}
\begin{proof} (1): \
By our definitions
$$M:=H_+(\om) \leftthreetimes_{\pi} G=((\R \times E \leftthreetimes_{\om^{\nabla}} F)
\leftthreetimes_{\a} \R_+) \leftthreetimes_{\pi} G.$$ Let $\g_1
\subset \g$ be a coarser Hausdorff group topology on $M$. Theorem
\ref{t:plus} establishes the minimality of $H_+(\om)=(\R \times E
\leftthreetimes_{\om} F) \leftthreetimes_{\a} \R_+$. In
particular, $\g_1$ and $\g$ agree on the subgroup $H:=H(\om)=\R
\times E \leftthreetimes_{\om^{\nabla}} F$.

By Lemma \ref{l:2.5} the natural projections $H \to E$ and $H \to
F$ are central. Since the birepresentation of $G$ in $\om$ is
t-exact we can apply Proposition \ref{p:2.7} (with $Y_1:=E, \
Y_2:=F$) to the group $H \leftthreetimes G$. It follows that
$\g_1$ agrees with $\g$ on the subgroup $H \leftthreetimes G$ of
$M$.


It is important now that $(M,\g)$ is an internal topological
semidirect product (see \cite[Section 6]{RD}) of $H
\leftthreetimes G$ with $\R_+$ (observe that $H \leftthreetimes G$
is a normal subgroup of $M=(H \leftthreetimes G) \cdot \R_+$ and
$(H \leftthreetimes G) \cap \R_+ = \{e_M\}$). This presentation
enables us to apply Proposition \ref{p:2.7}; this time in the
following situation: $G:=\R_+, \ X:=H \leftthreetimes G, \ Y:=\R$
and $q: X \to Y$ is the natural projection. Clearly $q: H
\leftthreetimes G \to \R$ is a continuous central retraction (with
respect to the topologies $\g|_{H \leftthreetimes G}$ and
$\g|_{\R}$).
The action of $\R_+$ on $\R$ is t-exact (Remark
\ref{r:definitionsand}.2). As we mentioned above $\g_1$ agrees
with $\g$ on the subgroup $X:=H \leftthreetimes G$ of $M$. So all
requirements of Proposition \ref{p:2.7} are fulfilled for the
topological group $(M,\g)$ and the coarser topology $\g_1 \subset
\g$. As a conclusion we get $\g_1=\g$.

(2): \ $ker(p)=H_+(\om)$ is minimal by Theorem \ref{t:plus}.
\end{proof}



\bibliographystyle{amsplain}

\end{document}